\documentclass[12pt,a4paper]{article}

\usepackage{amssymb,amsmath,amsthm,mathrsfs,cite}
\usepackage{color} 
\usepackage{verbatim}
\usepackage{soul}
\usepackage{multirow}
\usepackage{booktabs}
\usepackage{caption}
\usepackage{graphics}
\usepackage{threeparttable}
\usepackage{pifont}
\usepackage{float}
\usepackage{cases}
\usepackage{setspace}
\usepackage{ulem,cancel,bm}
\usepackage{makecell}
\soulregister{\cite}7 
\soulregister{\ref}7 
\soulregister{\bf}7
\soulregister{\it}7

\newtheorem{mythm}{Theorem}
\newtheorem{mylem}{Lemma}

\newtheorem{myrmk}{Remark}
\newtheorem{myprop}{Proposition}
\newtheorem{myexam}{Example}

\renewcommand{\arraystretch}{1.3}

\definecolor{gray}{gray}{0.55}

\begin{document}
	
	\begin{center}
		{\Large {\bf On Substochastic Inverse Eigenvalue Problems with the Corresponding Eigenvector Constraints}}\\[0.1in]
		Yujie Liu\footnote{Academy of Mathematics and Systems Science, Chinese Academy of Sciences, Beijing, 100190, P.R. China; School of Mathematical Sciences, University of Chinese Academy of Sciences, Beijing, 100049, China; yujieliu@amss.ac.cn} $\bullet$ Dacheng Yao\footnote{Academy of Mathematics and Systems Science, Chinese Academy of Sciences, Beijing, 100190, China; School of Mathematical Sciences, University of Chinese Academy of Sciences, Beijing, 100049, China; dachengyao@amss.ac.cn} $\bullet$ Hanqin Zhang\footnote{Department of Analytics \& Operations, National University of Singapore, 119245, Singapore; bizzhq@nus.edu.sg}
	\end{center}
	
	\noindent
	{\bf Abstract:} We consider the inverse eigenvalue problem of constructing a substochastic matrix from the given spectrum parameters with the corresponding eigenvector constraints. This substochastic inverse eigenvalue problem (SstIEP) with the specific eigenvector constraints is formulated into a nonconvex optimization problem (NcOP). The solvability for SstIEP with the specific eigenvector constraints is equivalent to identify the attainability of a zero optimal value for the formulated NcOP.  When the optimal objective value is zero, the corresponding optimal solution to the formulated NcOP is just the substochastic matrix desired to be constructed. We develop the alternating minimization algorithm to solve the formulated NcOP, and its convergence is established  by developing a novel method to obtain the boundedness of the optimal solution. Some numerical experiments are conducted to demonstrate the efficiency of the proposed method.
	
	\
	
	
	\
	             
	
	\section{Introduction}
	An $n$-by-$n$ nonnegative matrix $A=(a_{ij})_{n\times n}$ is said to be a stochastic matrix if all its row summations are one, and is called a substochastic matrix if all its row summations are less than or equal to one, but at least one of them is less than one. The stochastic inverse eigenvalue problem (StIEP) is to address the issue about whether there exists a stochastic matrix such that its spectrum is a given self-conjugate set of complex numbers $\{\lambda_1,\lambda_2,\ldots,\lambda_n\}$ (that is, solvability condition), and further if such a stochastic matrix exists, how to compute it. StIEP has many applications in applied probability, statistics, and engineering. In applied probability and statistics, for instance, StIEP is equivalent to construct a finite-state discrete-time Markov chain such that its transition probability  matrix is the corresponding constructed stochastic matrix. As any nonnegative matrix with positive maximal eigenvalue and positive maximal eigenvector can be transformed into a stochastic matrix by a diagonal similarity transformation (see \cite{chu-golub-book} or \cite{minc1988}), with the help of the well-studied nonnegative inverse eigenvalue problem and the Perron-Frobenius theorem, people can successfully identify some solvability conditions and develop algorithms to numerically compute the corresponding stochastic matrix, see, for instance, \cite{chu-guo-1998, yao-bai-2015, zhao-bai-2015, steidl2022}.
	
 In many real applications, however, we often need to solve the problem about whether there exists an $n$-by-$n$ substochastic matrix with some structural constraints on its corresponding eigenvectors such that not only is its spectrum a given self-conjugate set of complex numbers $\{\lambda_1,\lambda_2,\ldots,\lambda_n\}$ but its eigenvectors also satisfy the specified structural constraints. 
We call this problem the {\sf substochastic inverse eigenvalue problem (SstIEP) with the specified eigenvectors properties (SEP)}. 
For a given self-conjugate set of finite complex numbers, for example, in applied probability and engineering, people are interested in whether there exists an $(n+1)$-state discrete-time Markov chain with each state being transient except one of the states being an absorbing state such that the spectrum of its transition probability matrix among the $n$ transient states is the given self-conjugate set, at the same time, the corresponding eigenvectors satisfy some prior properties as well,  which is known as the {\it minimal phase-type distribution representation problem}, see \cite{he2014, he2008, pilungan2009}. Moreover, people need to find an effective way to compute it when such substochastic matrix  exists, see \cite{latouche1999}.
	
In this paper, we focus on SstIEP with SEP, and discuss its 
application to the minimal phase-type distribution representation problem. We specifically consider the scenario where the elements of the given set  $\{\lambda_1,\lambda_2,\ldots,\lambda_n\}$ are real and distinct. Actually, the existing research on the minimal phase-type distribution representation problem reveals that the location of the elements within the given set is crucial for identifying the corresponding Markov chain and even the simple case in which the elements are real and distinct is still challenging, see, e.g., \cite{commault2003, maier1991}.
	Inspired by this challenge, we consider such a problem, namely, for a given set of real and distinct numbers $\{\lambda_1,\lambda_2,\ldots,\lambda_n\}$, {\it when does there exist an $n$-by-$n$ substochastic matrix such that its spectrum is $\{\lambda_1,\lambda_2,\ldots,\lambda_n\}$, and its eigenvectors satisfy some specific properties}?  And when such substochastic matrix exists, {\it how do we effectively compute it}?
	
To address the SstIEP with SEP and the issue of effective computation, based on the framework of the alternating minimization ({\bf AM}) method from optimization theory,
(see \cite{bertsekas2015parallel, ye2021}), we propose a systematic method (algorithm) to simultaneously solve these two problems by formulating SstIEP with SEP into a nonconvex optimization problem with an unbounded convex feasible region.
	To make the {\bf AM} method work effectively for the nonconvex optimization problems, we need to prove the convergence of the proposed algorithm.
	The existing convergence analysis usually requires the compact level set or the quasi convexity of the objective function (see, e.g., \cite{beck2015, grippo2000convergence, hong2017iteration, tseng2001convergence}), while our nonconvex optimization problem does not satisfy any one of them. To obtain the convergence, we show the boundedness of the optimal solutions to this formulated nonconvex optimization problem with the unbounded convex feasible region. To establish the boundedness, we develop a new approach to estimate the lower bound of the determinant given by the matrix that characterizes the convex feasible region (CFR matrix for short). The approach consists of three steps: Doing perturbation on CFR matrix; Augmenting the perturbed CFR matrix while keeping the value of its corresponding determinant unchanged; and Coupling the submatrices from CFR matrix properly. Some numerical experiments are conducted to illustrate the effectiveness of the algorithm proposed here to solve SstIEP with SEP.
	
	Throughout this paper, ${\cal R}^n$ denotes the $n$-dimensional real space, and ${\cal R}_{+}^n$ is the set of  $n$-dimensional nonnegative real vectors. ${\cal R}^{n\times m}$ is the set of all $n$-by-$m$ real matrices, and
	${\cal R}_+^{n\times m}$ is the set of all nonnegative $n$-by-$m$ real matrices. The transpose of a vector or a matrix is denoted by appending a superscript ``$\top$".
	All the vectors in this paper are understood to be row vectors. For $a=(a_1,\ldots,a_n)\in{\cal R}^n$, its norm denoted by $\|a\|$ is defined by
	$\max_{1\leq i\leq n}\{|a_i|\}$, and for $A=(a_{ij})_{n\times n}\in {\cal R}^{n\times n}$, its norm written as $\|A\|$ is given by $\max_{1\leq i,j\leq n}\{|a_{ij}|\}$ while its  Frobenius norm is represented by  $\|A\|_F$. For $A\in {\cal R}^{n\times n}$, $A_{(i)}$ represents its $i$th row vector, $(a_{i1},\ldots,a_{in})$, with $1\leq i\leq n$,  $A^{(j)}$ represents its column vector, $(a_{1j},a_{2j},\cdots,a_{nj})$, with $1\leq j\leq n$,
	$\textrm{adj}(A)$ denotes its adjugate matrix, and $\det(A)$ is used to denote its determinant, and for $a=(a_1,\cdots,a_n)\in {\cal R}^n$, $\textrm{diag}(a)$ represents an $n$-by-$n$ diagonal matrix with the diagonal elements being $\{a_i: 1\leq i\leq n\}$.  Let $I$ be the identity matrix, ${\bf 1}$ the vector with each component being one, and ${\bf 0}$ the vector whose whole components are zero. 
	
	The paper is organized as follows:  Section \ref{formulation} formulates our problem and provides some preliminary results. In Section \ref{ADMM}, we develop the {\bf AM} method for our problem, and prove the convergence of the {\bf AM} method. Section \ref{numerical} is devoted to the numerical experiments, and the paper is concluded in Section \ref{conc}.
	
	\section{Problem Formulation}\label{formulation}
	Given $\lambda=(\lambda_1,\lambda_2,...,\lambda_n)\in {\cal R}^n$ and $\beta=(\beta_1,\beta_2,...,\beta_n)\in {\cal R}^n$ with $n\geq 2$, 
	\begin{align}
		& 1>\lambda_1>\lambda_2>\cdots>\lambda_n, \ \lambda_1>\max_{2\leq \ell\leq n}\{|\lambda_\ell|\};\label{1212-11}\\
		& \beta\cdot {\bf 1}^\top =1, \ \Pi_{\ell=1}^n\beta_\ell \neq0,\label{1212-12}
	\end{align}
	we have the following two problems closely related to the {\it substochastic inverse eigenvalue problem} and the {\it minimal phase-type distribution representation problem}. \\[-0.05in]
	
	\noindent
	{\bf Problem 1}: {\it Suppose that the prescribed numbers $\{\lambda_1,\ldots,\lambda_n\}$ given by {\rm (\ref{1212-11})}
is a realizable spectrum of some substochastic matrix. Do there exist $A\in {\cal R}^{n\times n}$ and $P\in {\cal R}^{n\times n}$ such that they satisfy the following constraints}?
	\begin{eqnarray}
		{\bf (P1):}\left\{
		\begin{array}{lllll}
			&PA= \textrm{diag}(\lambda)\cdot P;\\
			&P\cdot {\bf 1}^\top={\bf 1}^\top; \\
			&\beta\cdot   P\geq {\bf 0}; \\
			&A  \ \ \mbox{is nonnegative};\\
			&A \cdot {\bf 1}^\top \leq {\bf 1}^\top.
		\end{array}
		\right.
		\label{1212-13}
	\end{eqnarray}
	{\bf Problem 2}: {\it If there are $A$ and $P$ such that they satisfy {\rm (\ref{1212-13})}, then how can we numerically compute them}?

	\begin{myrmk}\label{rem-1212-1} 
		{\rm Given $n$ $(n\geq 2)$ distinct real numbers $\lambda_i$ $(1\leq i\leq n)$ with $\lambda_1$ being the largest,  $0<\lambda_1\leq 1$ is a necessary condition for $\{\lambda_1,\ldots,\lambda_n\}$
			to be the eigenvalues of an $n$-by-$n$ substochastic matrix from the Perron-Frobenius theorem (see Theorems 8.1.22 and 8.4.4 in \cite{horn2013}). However, the conditions in {\rm (\ref{1212-11})}, namely, $1>\lambda_1>\lambda_2>\cdots>\lambda_n$ and $\lambda_1>\max_{2\leq \ell\leq n}\{|\lambda_\ell|\}$, are neither necessary nor
			sufficient for $\{\lambda_1,\ldots,\lambda_n\}$
			to be the eigenvalues of an $n$-by-$n$ substochastic matrix.
			From $PA= \textrm{diag}(\lambda)P$ and $P\cdot {\bf 1}^\top={\bf 1}^\top$, we know that $P_{(\ell)}$ is the eigenvector of $A$ corresponding eigenvalue of $\lambda_\ell$, and is not orthogonal to ${\bf 1}$. It follows from {\rm (\ref{1212-11})} that
			$n$ vectors $P_{(\ell)}$ $(1\leq \ell\leq n)$ are independent. Thus, $P$ is nonsingular. Hence, $A=P^{-1}\cdot  \textrm{diag}(\lambda)\cdot P$.
			In view of $A \in {\cal R}^{n\times n}_+$ and $A \cdot {\bf 1}^\top \leq {\bf 1}^\top$ with {\rm (\ref{1212-11})}, this shows us that {\bf Problem 1 (P1)} is a
			{\it substochastic inverse eigenvalue problem  with structural constraints}: $PA=\textrm{diag}(\lambda)P$,  $P\cdot {\bf 1}^\top={\bf 1}^\top$, and $\beta\cdot P\geq {\bf 0}$, see \cite{chu-golub-book}. In other words, the structural constraints are related to the eigenvectors of the matrix determined by the substochastic inverse eigenvalue problem. 
			\hfill$\Box$}
	\end{myrmk}
	
	\begin{myrmk}\label{rem-1212-2}
		{\rm Following the paper's motivation discussed in the introduction, we now demonstrate how {\bf Problem 1} and {\bf Problem 2} are related to the minimal phase-type distribution representation problem, which is a long-lasting open problem in applied probability area. Given a rational function with real coefficients:
			\begin{align}\label{ori_f}
				& f(z)=\frac{p(z)}{q(z)}=\frac{p_nz^n+p_{n-1}z^{n-1}+\cdots+p_1z}{q_nz^n+q_{n-1}z^{n-1}+\cdots+q_1z+1},
			\end{align}
			where $p(z)$ and $q(z)$ are coprime, $p_n^2+q_n^2\neq 0$, and $p(1)=q(1)$,
			in applied probability area, the minimal phase-type distribution representation problem is whether there exists a discrete-time Markov chain with state space $\{0,1,\ldots,n\}$ such that
			(i) state $0$ is an absorbing state and the other states are transient; and  (ii) the moment generating function of its absorbing time is $f(z)$, see \cite{buchholz2014, he2014}.
			The transition probability matrix of the Markov chains with property (i) must have the following form:
			\begin{align*}
				\begin{bmatrix}
					1&{\bf 0}\\
					(I-A){\bf 1}^\top&A
				\end{bmatrix} \ \mbox{with $A\in {\cal R}_+^{n\times n}$ being a substochastic matrix}.
			\end{align*}	
			The initial distribution of such a Markov chain with property (ii) on states $\{1,\ldots,n\}$  is $\alpha\in {\cal R}^n_+$ with $\alpha\cdot{\bf 1}^{\top}= 1$.
			It follows from formula (2.14) of \cite{latouche1999} that the moment generating function of the absorbing time of the Markov chain with properties (i)-(ii) is
			$z\alpha(I-zA)^{-1}(I-A){\bf 1}^\top$.	
			
			Thus this minimal phase-type representation problem for the given rational function $f(z)$ in (\ref{ori_f}) is equivalent to verify whether there exist a substochastic matrix $A\in {\cal R}_+^{n\times n}$ and $\alpha\in {\cal R}^n_+$ with $\alpha\cdot{\bf 1}^{\top}= 1$ such that
			\begin{align}
				f(z)=z\alpha(I-zA)^{-1}(I-A){\bf 1}^\top.
				\label{1212-15}
			\end{align}
			Note that
			\[
			z\alpha(I-zA)^{-1}(I-A){\bf 1}^\top
			=\frac{z\alpha \cdot \textrm{adj}(I-zA)(I-A){\bf 1}^\top}{ \det(I-zA) }.
			\]
			It follows from (\ref{ori_f})-(\ref{1212-15}) that
			\begin{eqnarray}
				\frac{z\alpha\cdot  \textrm{adj}(I-zA)(I-A){\bf 1}^\top}{ \det(I-zA)}=\frac{p_nz^n+p_{n-1}z^{n-1}+\cdots+p_1z}{q_nz^n+q_{n-1}z^{n-1}+\cdots+q_1z+1}=\frac{p(z)}{q(z)}.
				\label{23-7-10}
			\end{eqnarray}
			The denominator and numerator on the left-hand-side of (\ref{23-7-10}) are two polynomials, both of their degrees are less than or equal to $n$. Since
			$\det(I-zA)$ equals to one when $z=0$, $p(z)$ and $q(z)$ are coprime, and $p_n^2+q_n^2\neq 0$,  thus (\ref{23-7-10}) implies
			\begin{eqnarray}
				 \det(I-zA)=q(z).\label{1212-14-july-1}
			\end{eqnarray}
			
			To see the relation between {\bf Problem 1} in (\ref{1212-13}) and the minimal phase-type representation, take a special case for $f(z)$ in \eqref{ori_f}, namely,  its denominator,   $q(z)$, has the form of $q(z)=(1-\lambda_1z)(1-\lambda_2z)\cdots(1-\lambda_nz)$ with $\{\lambda_1,\ldots,\lambda_n\}$ given by (\ref{1212-11}). Then  $f(z)$ given by \eqref{ori_f} can be, by the partial fraction decomposition, written as
			\begin{align}
				f(z)=z \Big(\frac{r_1}{1-\lambda_1z}+\frac{r_2}{1-\lambda_2z}+\cdots+ \frac{r_n}{1-\lambda_nz}\Big) \ \mbox{with} \ r_\ell\neq 0 \ \mbox{for $1\leq \ell\leq n$}. \label{1212-14}
			\end{align}
			
			First we demonstrate that if the minimal phase-type representation, $(\alpha, A)$,  exists such that {\rm (\ref{1212-15})} holds with $f(z)$ given by (\ref{1212-14}),
			then we can identify $\beta$ satisfying (\ref{1212-12}) and find $P$ such that $A$ and $P$ together make (\ref{1212-13}) hold with  $\beta$.
			To this end, by (\ref{1212-14-july-1}), the spectrum of $A$ is $\{\lambda_1,\lambda_2,\dots,\lambda_n\}$. Thus, there exists a nonsingular $P\in {\cal R}^{n\times n}$ such that
			$P A= \textrm{diag}(\lambda)P$.
			Then we have $A=P^{-1} \textrm{diag}(\lambda)P$. This together with {\rm (\ref{1212-15}) and (\ref{1212-14})} gives us
			\begin{align}
				& z \Big(\frac{r_1}{1-\lambda_1z}+\frac{r_2}{1-\lambda_2z}+\cdots+ \frac{r_n}{1-\lambda_nz}\Big)\nonumber\\
				& \ \ \ =z\alpha \Big(I-z P^{-1} \textrm{diag}(\lambda)P\Big)^{-1}\Big(I-P^{-1} \textrm{diag}(\lambda)P\Big){\bf 1}^\top\label{1212-17} \\
				& \ \ \ =z\alpha P^{-1}\Big(I-z  \textrm{diag}(\lambda)\Big)^{-1}\Big(I- \textrm{diag}(\lambda)\Big)P \cdot {\bf 1}^\top. \nonumber
			\end{align}
			Note that
			\begin{align}\label{diag}
				\Big(I-z \textrm{diag}(\lambda)\Big)^{-1}\Big(I- \textrm{diag}(\lambda)\Big )= \textrm{diag}\Big(\frac{1-\lambda_1}{1-z\lambda_1},\ldots,\frac{1-\lambda_n}{1-z\lambda_n}\Big),
			\end{align}
			which together with (\ref{1212-14})-(\ref{1212-17}) gives us each component of $P{\bf 1}^\top$ is not zero.
			Since each row of $P$ is an eigenvector of $A$, then without loss of generality, we can assume
			$P{\bf 1}^\top={\bf 1}^\top$,
	 and again use (\ref{1212-17})-(\ref{diag}),
		\begin{align*}
			z \Big(\frac{r_1}{1-\lambda_1z}+\frac{r_2}{1-\lambda_2z}+\cdots+ \frac{r_n}{1-\lambda_nz}\Big)
            & =z\alpha P^{-1} \textrm{diag}\Big(\frac{1-\lambda_1}{1-z\lambda_1},\ldots,\frac{1-\lambda_n}{1-z\lambda_n}\Big)P \cdot {\bf 1}^\top\\
            & =z\alpha P^{-1} \textrm{diag}\Big(\frac{1-\lambda_1}{1-z\lambda_1},\ldots,\frac{1-\lambda_n}{1-z\lambda_n}\Big){\bf 1}^\top.
		\end{align*}
Noting that $\alpha P^{-1}$ has nothing to do with $z$, and comparing two rational functions of $z$ given by the two sides of the above equation, we have
$\alpha P^{-1}=(\frac{r_1}{1-\lambda_1},\frac{r_2}{1-\lambda_2},\ldots,\frac{r_n}{1-\lambda_n})$.
		Let	
		\begin{align}\label{beta}
			\beta=\bigg(\frac{r_1}{1-\lambda_1},\frac{r_2}{1-\lambda_2},\ldots,\frac{r_n}{1-\lambda_n}\bigg).
		\end{align}
 Then \eqref{1212-12} holds from \eqref{1212-14} and $f(1)=1$, and the expression of $\alpha P^{-1}$  and (\ref{beta}) give that $\beta=\alpha P^{-1}$. It implies that $\beta P=\alpha\geq {\bf 0}$. Hence we know that if such a minimal phase-type representation problem has a solution, then there exist  $A\in {\cal R}^{n\times n}$ and $P\in {\cal R}^{n\times n}$ such that \eqref{1212-13} holds with $\beta$ given by \eqref{beta}.
			
			Next we show that for $\lambda$ and $\beta$ given by \eqref{beta}, if there exist  $A\in {\cal R}^{n\times n}$ and $P\in {\cal R}^{n\times n}$ such that \eqref{1212-13} holds, then
			$(\alpha, A)$ with $\alpha=\beta\cdot P$ is the minimal phase-type representation for $f(z)$ given by \eqref{1212-14}, that is, we find that the pair $(\alpha, A)$ satisfies (\ref{1212-15}) with $f(z)$ given by \eqref{1212-14}. To get this, note that $\beta$ given by \eqref{beta} satisfies \eqref{1212-12} by the fact that $f(1)=1$ and \eqref{1212-14}.
			We observe that $\beta \cdot P\geq {\bf 0}$ together with $P\cdot {\bf 1}^\top={\bf 1}^\top$ and \eqref{1212-12} gives that $\alpha\in {\cal R}^n_+$, and  $\alpha\cdot{\bf 1}^{\top}= 1$.
			Remark \ref{rem-1212-1} gives that $P$ is nonsingular and $A=P^{-1}\cdot \textrm{diag}(\lambda)\cdot P$.
			Furthermore, from \eqref{diag}-\eqref{beta} and $P\cdot {\bf 1}^\top={\bf 1}^\top$, we have
			\begin{align*}
				& z \Big(\frac{r_1}{1-\lambda_1z}+\frac{r_2}{1-\lambda_2z}+\cdots+ \frac{r_n}{1-\lambda_nz}\Big)\\
				& \ \ \ =z\beta\Big(I-z  \textrm{diag}(\lambda)\Big)^{-1}\Big(I- \textrm{diag}(\lambda)\Big)\cdot {\bf 1}^\top\\
				& \ \ \ =z\beta\Big(I-z  \textrm{diag}(\lambda)\Big)^{-1}\Big(I- \textrm{diag}(\lambda)\Big)P \cdot {\bf 1}^\top \ \mbox{(by $P\cdot {\bf 1}^\top={\bf 1}^\top$)}\\
				& \ \ \ =z\beta P\Big(I-z P^{-1} \textrm{diag}(\lambda)P\Big)^{-1}\Big(I-P^{-1} \textrm{diag}(\lambda)P\Big) \cdot {\bf 1}^\top\\
				& \ \ \ =z\alpha(I-z A)^{-1}(I-A) {\bf 1}^\top.
			\end{align*}
			Here the last equation follows from  $\alpha=\beta\cdot P$ and $A=P^{-1}\cdot  \textrm{diag}(\lambda)\cdot P$. This implies that the minimal phase-type representation problem for $f(z)$ given by (\ref{1212-14}) has a solution.
			
			In sum, we know that the minimal phase-type representation problem for $f(z)$ given by (\ref{1212-14}) is equivalent
			to {\bf Problem 1} and {\bf Problem 2}.} \hfill$\Box$
	\end{myrmk}
	
	In the following, we first formulate {\bf Problem 1} as a nonconvex optimization problem,  then use the {\bf AM} approach in nonconvex optimization theory to simultaneously solve {\bf Problem 1} and {\bf Problem 2}. For the {\bf AM} approach, see \cite{bertsekas2015parallel}.
	
	\section{ {\bf AM} Approach}\label{ADMM}
	First we reformulate {\bf Problem 1 (P1)} as an optimization problem:
	\begin{eqnarray}
		{\bf (OP):}\left\{
		\begin{array}{ll}
			& \min_{P,A} \Big\|PA-\textrm{diag}(\lambda)P \Big\|^2_{F} \label{july-15-1}\\
			& {\sf s.t.} \ \left\{
			\begin{array}{llll}
				& P {\bf 1}^\top={\bf 1}^\top;
				\\
				&\beta  P\geq {\bf 0};\\
				& A \ \mbox{is nonnegative};\\
				&A {\bf 1}^\top \leq {\bf 1}^\top.
			\end{array}
			\right.
		\end{array}
		\right. \label{july-15-2}
	\end{eqnarray}
	We know that {\bf Problem 1} given by (\ref{1212-13}) is solvable if and only if the above optimization problem ({\bf OP}) given by (\ref{july-15-2}) has a solution such that the optimal objective value is zero. As mentioned in the above section, we use the {\bf AM} approach to explore the optimal solution of the problem {\bf (OP)}.\\[0.15in]
	{\bf Step-a:} At the beginning, we choose $A_0\in {\cal R}_+^{n\times n}$ with $A_0{\bf 1}^\top\leq {\bf 1}^\top$, and solve the following optimization problem  with $A=A_0$:
	\begin{eqnarray}
		{\bf (OP[A])}: \left \{
		\begin{array}{ll}
			&\min_{P} \Big\|P A- \textrm{diag}(\lambda)P \Big\|^2_{F}\\
			& {\sf s.t.} \ \left\{
			\begin{array}{ll}
				& P {\bf 1}^\top ={\bf 1}^\top;  \\
				&\beta P\geq {\bf 0}.
			\end{array}
			\right.
		\end{array}
		\right. \label{july-15-4}
	\end{eqnarray}
	
	\noindent
	{\bf Step-b:} Let $P_0$ denote the solution of the optimization problem {\bf (OP[A])} given by (\ref{july-15-4}). Then we solve another optimization problem with $P=P_0$:
	\begin{eqnarray}
		{\bf (OP[P])}: \left \{
		\begin{array}{ll}
			&\min_{A} \Big\|P A-\textrm{diag}(\lambda) P \Big\|^2_{F} \label{july-15-5}\\
			& {\sf s.t.} \ \left\{
			\begin{array}{ll}
				& A \ \mbox{is nonnegative}; \\
				& A {\bf 1}^\top \leq {\bf 1}^\top.
			\end{array}
			\right.
		\end{array}
		\right. \label{july-15-6}
	\end{eqnarray}
	
	\noindent
	{\bf Step-c:} Let $A_1$ be the solution to (\ref{july-15-6}). Then solve the optimization problem {\bf (OP[A])} given by  (\ref{july-15-4}) with replacing  $A$ by $A_1$, and its solution is denoted by $P_1$. Again, we solve the optimization problem {\bf (OP[P])} given by  (\ref{july-15-6}) with replacing $P$ by $P_1$, and denote its solution by $A_2$. Repeating {\bf Step-a} and {\bf Step-b}, we generate a sequence $\{(A_\ell,P_\ell): \ell\geq 0\}$. \\[-0.05in]
	
	To make the above steps workable, we establish the following lemmas.
	
	\begin{mylem} \label{objective-convexity}
		For any given $A=(a_{ij})_{n\times n}\in {\cal R}^{n\times n}_+$ with $A{\bf 1}^\top\leq {\bf 1}^\top$,
		{\rm $\Big\|PA-\textrm{diag}(\lambda)P \Big\|^2_{F}$ }
		is convex with respect to $n^2$ variables $\{p_{ij}: 1\leq i,j\leq n\}$ given by the entries of matrix $P=(p_{ij})_{n\times n}$.
		Moreover, {\rm $\Big\|PA- \textrm{diag}(\lambda)P \Big\|^2_{F}$} can be written as
		\begin{align}\label{ob_exp}
			\Big(P_{(1)},P_{(2)},...,P_{(n)} \Big)
			\begin{bmatrix}
				A[\lambda_1]& & &\\
				& A[\lambda_2] & &\\
				& & \ddots & \\
				& & &A[\lambda_n]
			\end{bmatrix}
			\left(\begin{array}{cccc}
				P_{(1)}^{\top}\\
				P_{(2)}^{\top}\\
				\vdots\\
				P_{(n)}^{\top}
			\end{array}
			\right),
		\end{align}
		where
		\begin{align}
			A[x]=(A-xI)(A-xI)^{\top} \ \mbox{for any $x\in {\cal R}$}.
			\label{1212-50}
		\end{align}
		Similarly, for any given {\rm $P=(p_{ij})_{n\times n}\in {\cal R}^{n\times n}$, $\Big\|PA- \textrm{diag}(\lambda)P \Big\|^2_{F}$ }
		is convex with respect to $n^2$ variables $\{a_{ij}: 1\leq i,j\leq n\}$ given by the entries of matrix $A=(a_{ij})_{n\times n}$.
		Further, {\rm $\Big\|PA- \textrm{diag}(\lambda)P \Big\|^2_{F}$}  can be equivalently written into another expression
		\begin{align}\label{obja}
			\sum_{i,j=1}^n (A^{(j)})^\top P_{(i)}^\top P_{(i)}  A^{(j)}-2\lambda_i p_{ij} P_{(i)} A^{(j)}+\lambda_i^2 p_{ij}^2,
		\end{align}
		where $A^{(j)}=(a_{1j},a_{2j},...,a_{nj})^{\top}$.
	\end{mylem}
	
	\noindent
	{\it Proof}: These two convexities are straightforward. Now prove $\Big\|PA- \textrm{diag}(\lambda)P \Big\|^2_{F}$ can be written into the above two different expressions.
	Note that
	\begin{align*}
		PA- \textrm{diag}(\lambda)P=\left(\begin{array}{cccc}
			P_{(1)}\\
			P_{(2)}\\
			\vdots\\
			P_{(n)}
		\end{array}
		\right)\cdot A-\left(\begin{array}{cccc}
			\lambda_1P_{(1)}\\
			\lambda_2P_{(2)}\\
			\vdots\\
			\lambda_nP_{(n)}
		\end{array}
		\right)=\left(\begin{array}{cccc}
			P_{(1)}(A-\lambda_1I)\\
			P_{(2)}(A-\lambda_2I)\\
			\vdots\\
			P_{(n)}(A-\lambda_nI)
		\end{array}
		\right).
	\end{align*}
	Thus
	\begin{align*}
		\Big\|PA-\textrm{diag}(\lambda)P \Big\|^2_{F}=\sum_{i=1}^nP_{(i)}(A-\lambda_iI)(A-\lambda_iI)^{\top}P_{(i)}^{\top}.
	\end{align*}
	This gives (\ref{ob_exp}). The second expression directly follows from that for any $1\leq i,j\leq n$,
	\begin{eqnarray*}
		\Big((PA- \textrm{diag}(\lambda)P)_{ij}\Big)^2&=&\Big(P_{(i)}A^{(j)}-\lambda_ip_{ij}\Big)^2\nonumber\\
		&=& (A^{(j)})^\top P_{(i)}^\top P_{(i)}  A^{(j)}-2\lambda_i p_{ij} P_{(i)} A^{(j)}+\lambda_i^2 p_{ij}^2.
	\end{eqnarray*}
	Hence, we have the lemma.
	\hfill$\Box$
	
	\begin{mylem}\label{Alter-optimal-solution}
		For any given $A\in {\cal R}_+^{n\times n}$ with $A{\bf 1}^\top\leq {\bf 1}^\top$, the convex optimization problem ${\bf (OP[A])}$ given by {\rm (\ref{july-15-4})} has an optimal solution$;$ For any given $P\in {\cal R}^{n\times n}$, the convex optimization problem ${\bf (OP[P])}$ given by {\rm (\ref{july-15-6})} has an optimal solution.
	\end{mylem}	
	
	\noindent
	{\it Proof:} We first prove the first part using Corollary 2.1 of \cite{lee2005quadratic}.
	To transform our convex optimization problem ${\bf (OP[A])}$ given by (\ref{july-15-4}) into the problem addressed in Corollary 2.1 of \cite{lee2005quadratic}, according to \eqref{ob_exp}, take $D$, $c$, $b$ and $A$ in Problem (2.1) on p.29 studied by Corollary 2.1 of \cite{lee2005quadratic} to be
	\begin{eqnarray*}
		D=\begin{bmatrix}
			A[\lambda_1]& & &\\
			& A[\lambda_2] & &\\
			& & \ddots & \\
			& & &A[\lambda_n]
		\end{bmatrix}; \ c={\bf 0}; \  b=\left(
		\begin{array}{ccc}
			{\bf 1}^\top\\
			-{\bf 1}^\top\\
			{\bf 0}^\top
		\end{array}
		\right); \ A=\begin{bmatrix}
			E\\
			-E\\
			I[\beta]
		\end{bmatrix}
	\end{eqnarray*}
	with
	\[
	E=\begin{bmatrix}
		{\bf 1} & & &\\
		&{\bf 1}& &\\
		&  & \ddots &\\
		& &  &{\bf 1}
	\end{bmatrix}\in {\cal R}^{n\times n^2}_+
	\ \mbox{and} \ \ I[\beta]=\begin{bmatrix}\beta_1I,\cdots,\beta_nI
	\end{bmatrix}\in {\cal R}^{n\times n^2}.
	\]
	Then by  Corollary 2.1 of \cite{lee2005quadratic}, it suffices to prove that (i) the constraint set of ${\bf (OP[A])}$ given by {\rm (\ref{july-15-4})} is nonempty; and (ii) if $V=(v_{ij})\in {\cal R}^{n\times n}$ satisfies $V {\bf 1}^\top ={\bf 0}^\top$, $\beta V\geq {\bf 0}$, and
	\begin{align}
		\Big(V_{(1)},V_{(2)},...,V_{(n)} \Big)
		\begin{bmatrix}
			A[\lambda_1]& & &\\
			& A[\lambda_2] & &\\
			& & \ddots & \\
			& & &A[\lambda_n]
		\end{bmatrix}
		\left(\begin{array}{cccc}
			V_{(1)}^{\top}\\
			V_{(2)}^{\top}\\
			\vdots\\
			V_{(n)}^{\top}
		\end{array}
		\right)=0,
		\label{1212-20}
	\end{align}
	then for $P$ satisfying $P {\bf 1}^\top ={\bf 1}^\top$ and $\beta P\geq {\bf 0}$, we have
	\begin{align}\label{1212-21}
		\Big(P_{(1)},P_{(2)},...,P_{(n)} \Big)
		\begin{bmatrix}
			A[\lambda_1]& & &\\
			& A[\lambda_2] & &\\
			& & \ddots & \\
			& & &A[\lambda_n]
		\end{bmatrix}
		\left(\begin{array}{cccc}
			V_{(1)}^{\top}\\
			V_{(2)}^{\top}\\
			\vdots\\
			V_{(n)}^{\top}
		\end{array}
		\right)\geq0.
	\end{align}
	It is direct to verify
	\begin{align}\label{delta-exp}
		\Delta=\begin{bmatrix}
			 \sum_{i=1}^n \beta_i/\beta_1 & -\sum_{i=2}^n \beta_i/\beta_1 & 0 & 0 & \cdots & 0 \\
			& \sum_{i=2}^n \beta_i/\beta_2  & -\sum_{i=3}^n\beta_i/ \beta_2 & 0 & \cdots & 0\\
			&        &  \sum_{i=3}^n \beta_i/\beta_3 & -\sum_{i=4}^n \beta_i/\beta_3 & \cdots & \vdots\\
			& & & \ddots & \ddots & 0\\
			& & & &\ddots & -\beta_n/\beta_{n-1}\\
			& & & & & 1
		\end{bmatrix}
	\end{align}
	satisfies the constraint set of ${\bf (OP[A])}$ given by {\rm (\ref{july-15-4})}, thus we have (i).
	
	Now we prove (ii).
	Note that the left-hand side of (\ref{1212-20}) can be written as
	$\sum_{\ell=1}^n V_{(\ell)} \Big(A-\lambda_\ell I\Big)\Big(A-\lambda_\ell I\Big)^{\top}V_{(\ell)}^\top$.
	It follows from \eqref{1212-20} that
	\[
	\Big(A-\lambda_\ell I\Big)^{\top}V_{(\ell)}^\top={\bf 0} \ \ \mbox{for each} \ \ell=1,\ldots,n.
	\]
	By the definition of $A[\lambda_i]$ given by (\ref{1212-50}), this implies
	\begin{align*}
		\begin{bmatrix}
			A[\lambda_1]& & &\\
			& A[\lambda_2] & &\\
			& & \ddots & \\
			& & &A[\lambda_n]
		\end{bmatrix}
		\left(\begin{array}{cccc}
			V_{(1)}^{\top}\\
			V_{(2)}^{\top}\\
			\vdots\\
			V_{(n)}^{\top}
		\end{array}
		\right)=0.
	\end{align*}	
	Then (\ref{1212-21}) holds, that is, we have (ii). Hence, the first part of the lemma is proved.
	
	For the second part, the conclusion is immediate because the constraint set given by $A\in{\cal R}^{n\times n}_+$ and $A{\bf 1}\leq {\bf 1}$  is compact and the objective function,
	$\|PA- \textrm{diag}(\lambda)P \|^2_{F}$, is continuous with respect to $n^2$ variables given by entries of $A$ for any fixed $P$. \hfill$\Box$\\[-0.05in]

	Lemma \ref{objective-convexity} shows that each of {\bf Step-a} and {\bf Step-b} is a convex optimization problem.
	Further, Lemma \ref{Alter-optimal-solution} makes that each step is well posed.
	To establish the workability of the
	{\bf AM} approach for the nonconvex optimization problem ({\bf OP}) given by (\ref{july-15-2}),
	we need to prove the convergence of $\{(A_\ell,P_\ell): \ell\geq 0\}$ and that the corresponding limit satisfies its Karush-Kuhn-Tucker ({\bf KKT}) conditions. To this end, define
	\[
	{\cal L}(A,P)=\Big\|PA-\textrm{diag}(\lambda)P \Big\|^2_{F}.
	\]
	
	We first look at the {\bf KKT} conditions corresponding to the convex optimization problem {\bf (OP[A])} given by (\ref{july-15-4})  and the convex optimization problem {\bf (OP[P])} given by (\ref{july-15-6}). They can be directly written as:
	\begin{align}
		{\bf KKT}(\bf OP[A]):	\left\{
		\begin{array}{lllllll}
			&\frac{\partial}{\partial p_{ij}} {\cal L}(A,P)+\alpha_{i} - \widetilde\alpha_{j}\cdot \beta_i=0, \ 1\leq i,j\leq n;\\
			& P{\bf 1}^\top={\bf 1}^\top;\\
			& \beta P\geq {\bf 0};\\
			& \widetilde\alpha_j \cdot \Big(\sum_{\ell=1}^n \beta_\ell p_{\ell j}\Big)=0, \ 1\leq j\leq n;\\
			& (\widetilde\alpha_1,\ldots, \widetilde\alpha_n)\in {\cal R}^n_+,
		\end{array}
		\right.
		\label{1212-40}
	\end{align}
	\begin{align}
		{\bf KKT (OP[P])}:	\left\{
		\begin{array}{lllllll}
			&\frac{\partial}{\partial a_{ij}} {\cal L}(A,P)+\pi_{i} -\widetilde \pi_{ij}=0, \ 1\leq i,j\leq n;\\
			& A{\bf 1}^\top\leq {\bf 1}^\top;\\
			& \pi_i \times \Big(\sum_{\ell=1}^n a_{i\ell}-1\Big)=0, \ 1\leq i\leq n;\\
			& A \ \mbox{is nonnegative};\\
			& \widetilde \pi_{ij}a_{ij}=0, \ 1\leq i,j\leq n;\\
			&(\pi_1,\ldots,\pi_n)\in {\cal R}^n_+ \ \mbox{and} \ (\widetilde{\pi}_{ij})_{n\times n}\in {\cal R}^{n\times n}_+.
		\end{array}
		\right.
		\label{1212-31}
	\end{align}
	Then the {\bf KKT} conditions corresponding to the nonconvex optimization problem $({\bf OP})$ given by (\ref{july-15-2}), denoted by {\bf KKT}($\bf OP$), are given by (\ref{1212-40})-(\ref{1212-31}).
	To obtain the convergence of $\{(A_\ell,P_\ell): \ell\geq 0\}$ given by {\bf Step-c}, we need to establish the boundedness of the optimal solution to {\bf (OP[A])} given by (\ref{july-15-4}) with any given $A \in {\cal R}^{n\times n}_+$ and $A{\bf 1}^\top\leq {\bf 1}^\top$. That is, we show there exists ${\rho}$ such that for any $A \in {\cal R}^{n\times n}_+$ with $A{\bf 1}^\top\leq {\bf 1}^\top$, the optimal solution to {\bf (OP[A])} given by (\ref{july-15-4}), denoted by $P_{\sf sln}(A)$, must satisfy
	\begin{align}
		\|P_{\sf sln}(A)\|\leq \rho.\label{uniform_bound}
	\end{align}
	
	Let $P=\Delta$ given by \eqref{delta-exp}. Note that $P$
	is a feasible solution to the optimization problem {\bf (OP[A])} given by (\ref{july-15-4}). Let  $\eta=\max\{|\beta_{\ell}|/|\beta_i|, {\ell}\geq i\}$. We have $|\sum_{{\ell}=i+1}^n \beta_{\ell}/\beta_i| \leq (n-i)\eta$ for $1\leq  \ell\leq n-1$. Moreover, by (\ref{1212-12}),
	\begin{align*}
		\|\Delta \|^2_{F}&=\sum_{i=1}^n \Big(\frac{\sum_{\ell=i}^n\beta_\ell}{\beta_{i}}\Big)^2+\sum_{i=2}^{n}
		\Big(\frac{\sum_{\ell=i}^n\beta_\ell}{\beta_{i-1}}\Big)^2\\
		&\leq \sum_{i=1}^n ((n-i)\eta+1)^2+\sum_{i=1}^{n-1}((n-i)\eta)^2\\
		&=\frac{(n-1)n(2n-1)}{3} \eta^2+n(n-1)\eta+n.
	\end{align*}
	This implies that for each $A \in {\cal R}^{n\times n}_+$ and $A{\bf 1}^\top\leq {\bf 1}^\top$,
	\begin{align*}
		{\cal L}(A, \Delta)&=\|\Delta\cdot A-\textrm{diag}(\lambda) \Delta\|_{F}^2\\
		&\leq\big( \| \Delta \cdot A \|_{F} +\| \textrm{diag}(\lambda) \Delta  \|_{F} \big)^2 \\
		&\leq\big( \| \Delta\|_{F} \cdot \|A \|_{F} + \lambda_1\| \Delta  \|_{F} \big)^2 \\
		&\leq (n+\lambda_1)^2\cdot \|\Delta \|^2_{F}\\
		&\leq (n+\lambda_1)^2\Big(\frac{(n-1)n(2n-1)}{3} \eta^2+n(n-1)\eta+n\Big).
	\end{align*}
	Denote
	\begin{align}
		\bar {\rho}=(n+\lambda_1)^2\Big(\frac{(n-1)n(2n-1)}{3} \eta^2+n(n-1)\eta+n\Big).
		\label{1212-37}
	\end{align}
	Then the optimal solution to {\bf (OP[A])} with $A \in {\cal R}^{n\times n}_+$ and $A{\bf 1}^\top\leq {\bf 1}^\top$ given by (\ref{july-15-4}) must satisfy
	\begin{align*}
		\|PA- \textrm{diag}(\lambda) P \|^2_ {F}=\sum_{i=1}^nP_{(i)}(A-\lambda_iI)(A-\lambda_iI)^{\top}P_{(i)}^{\top} \leq \bar {\rho},
	\end{align*}
	which gives that for $1\leq i \leq n$,
	
	\begin{align*}
		P_{(i)}(A-\lambda_i I)(A-\lambda_iI)^{\top}P_{(i)}^{\top} \leq \bar {\rho}.
	\end{align*}
	Then we have for $1\leq i \leq n$,
	\begin{align}\label{lu_xi}
		-\sqrt{\bar {\rho}} {\bf 1}  \leq P_{(i)}(A-\lambda_i I)\leq  \sqrt{\bar {\rho}} {\bf 1}.
	\end{align}
	In addition, $\beta P{\bf 1}^\top=\beta {\bf 1}^\top=1$ and $\beta P\geq \bm{0}$ imply that $\beta P\leq {\bf 1}$. This together with \eqref{lu_xi} implies that the optimal solutions to  {\bf (OP[A])}  with $A \in {\cal R}^{n\times n}_+$ and $A{\bf 1}^\top\leq {\bf 1}^\top$ given by (\ref{july-15-4}) must satisfy
	\begin{align}\label{add_con}
		\left(\begin{array}{ccccc}
			-\sqrt{\bar {\rho}}{\bf 1}^\top\\
			-\sqrt{\bar {\rho}}{\bf 1}^\top\\
			\vdots\\
			-\sqrt{\bar {\rho}}{\bf 1}^\top\\
			{\bf 0}^\top
		\end{array}\right)\leq
		\begin{bmatrix}
			A^{\top}-\lambda_1I & & &\\
			& A^{\top}-\lambda_2I & &\\
			& & \ddots &\\
			& & & A^{\top}-\lambda_nI\\
			\beta_1 I& \beta_2I &\cdots & \beta_n I
		\end{bmatrix}
		\left(\begin{array}{cccc}
			P_{(1)}^{\top}\\
			P_{(2)}^{\top}\\
			\vdots\\
			P_{(n)}^{\top}
		\end{array}\right)
		\leq \left(\begin{array}{ccccc}
			\sqrt{\bar {\rho}}{\bf 1}^\top\\
			\sqrt{\bar {\rho}}{\bf 1}^\top\\
			\vdots\\
			\sqrt{\bar {\rho}}{\bf 1}^\top\\
			{\bf 1}^\top
		\end{array}\right).
	\end{align}
	In order to get  \eqref{uniform_bound} from (\ref{add_con}), consider matrix $B^\top$ given by
	\begin{eqnarray}
		\begin{bmatrix}
			A^{\top}-\lambda_1I & & &\\
			& A^{\top}-\lambda_2I & &\\
			& & \ddots &\\
			& & & A^{\top}-\lambda_nI\\
			\beta_1 I& \beta_2I &\cdots & \beta_n I
		\end{bmatrix}. \label{12-10-1}
	\end{eqnarray}
	In view of (\ref{add_con}), we need to consider the singularity of $B^\top$. To simplify some formula expressions in our analysis, for any positive integer number $\ell$, we now define
	\begin{align*}
		[\ell]_n=\left\{
		\begin{array}{ll}
			\ell({\sf mod \ }n), & \mbox{if} \ \ell({\sf mod \ }n)\neq 0,\\
			n, & \mbox{if} \ \ell({\sf mod \ }n)=0.
		\end{array}
		\right.
	\end{align*}
	Let
	\begin{align}
		\gamma_i^{(1)}&=\beta^2_i, \ \ 1\leq i\leq n;\label{1212-51}\\
		\gamma_{i}^{(k+1)}&=\frac{\gamma_i^{(k)}\gamma_{[i+1]_n}^{(k)}(\lambda_i-\lambda_{[i+k]_n})^2}{\gamma_i^{(k)}+\gamma_{[i+1]_n}^{(k)}}, \ \ 1 \leq i\leq n \ \mbox{and} \ \ 1\leq k\leq n-2.\label{1212-52}
	\end{align}
	Since $\lambda_i\neq \lambda_{[i+k]_n}$ holds by \eqref{1212-11}, from \eqref{1212-12}, the $\gamma_{i}^{(k)}$ all stay positive for $1 \leq i\leq n$ and $1\leq k\leq n-1$.
	The non-singularity of $B^\top$ is obtained by establishing a lower bound on the determinant of $BB^{\top}$.
	\begin{myprop}\label{tech-prop} Consider $B^\top$ given by {\rm (\ref{12-10-1})}. We have
		\begin{align}\label{lb_B}
			\det(BB^{\top})\geq  2^{-n(n-2)}\bigg(\frac{\sum_{i\neq j}\gamma_{i}^{(n-1)}\gamma_{j}^{(n-1)}(\lambda_{[n+i-1]_n}-\lambda_{[n+j-1]_n})^2}{\sum_{i=1}^n \gamma_{i}^{(n-1)}}\bigg)^n
		\end{align}
		for any $A\in {\cal R}^{n\times n}_+$ with $A{\bf 1}^\top \leq {\bf 1}^\top$.
	\end{myprop}
	
	First we outline its proof. For a given matrix, usually it is hard to estimate a lower bound of its determinant.
	With a consideration of the possible singularity of some among $n$ matrices $A^\top-\lambda_iI$ $(1\leq i\leq n)$, we introduce perturbations for $A^\top-\lambda_iI$. Namely,
	for $t\in {\cal R}$, let
	\begin{align*}
		C(t) &= \begin{bmatrix}
			A[\lambda_1+t]+\beta_1^2I &\beta_1\beta_2I &\cdots &  \beta_1\beta_n  I\\
			\beta_1 \beta_2 I &A[\lambda_2+t]+\beta_2^2I &\cdots & \beta_2 \beta_nI\\
			\vdots& & \ddots &  \vdots\\
			\beta_1\beta_n  I&\beta_2 \beta_nI & \cdots&A[\lambda_n+t]+\beta_n^2I
		\end{bmatrix},
	\end{align*}
	where $A[\lambda_i+t]$ 
		is defined in (\ref{1212-50}).
	By the continuity of $ \det\big(C(t)\big)$ in $t$, we only need to consider $t$ such that $A[\lambda_i+t]$ is nonsingular. We properly augment one row and one column of $C(t)$ while keeping the determinant unchanged. Based on the augmented matrix of $C(t)$, using the Schur complement formula, $ \det\big(C(t)\big)$ can be written as the summation of products of $(n-1)$ terms of
	$A[\lambda_i+t]$ $(1\leq i\leq n)$. Among these $n$ products of $(n-1)$ terms of $A[\lambda_i+t]$ $(1\leq i\leq n)$, we do coupling for each two products and use the Minkowski inequality, then we obtain a lower bound for $ \det\big(C(t)\big)$ given by the summations of products of $(n-2)$ terms of $A[\lambda_i+t]$ $(1\leq i\leq n)$.  Again do coupling and by the Minkowski inequality, the lower bound of $ \det\big(C(t)\big)$ can be further reduced to the summations of products of $(n-3)$ terms of $A[\lambda_i+t]$ $(1\leq i\leq n)$. Repeat this and finally get the proposition.\\[-0.05in]
	
	\noindent
	{\it Proof} {\bf of Proposition \ref{tech-prop}:} \  First we have
	\begin{align*}
		 \det\big(C(t)\big) = \det\left(\begin{array}{ccccc}
			A[\lambda_1+t]+\beta_1^2I &\beta_1\beta_2I &\cdots &  \beta_1\beta_n  I& \beta_1 I\\
			\beta_1 \beta_2 I &A[\lambda_2+t]+\beta_2^2I &\cdots & \beta_2 \beta_nI & \beta_2 I\\
			\vdots& & \ddots &  \vdots& \vdots\\
			\beta_1\beta_n  I&\beta_2 \beta_nI & \cdots&A[\lambda_n+t]+\beta_n^2I& \beta_n I\\
			0 & 0 & \cdots& 0& I
		\end{array}
		 \right),
	\end{align*}
	and
	\begin{align}
		 \det(BB^\top)=\det\big(C(0)\big). \label{1212-9}
	\end{align}
	Moreover, we have
	\begin{align*}
		\det\big(C(t)\big) = \det\left(\begin{array}{ccccc}
			A[\lambda_1+t] &0 &\cdots &  0& \beta_1 I\\
			0  &A[\lambda_2+t] &\cdots & 0& \beta_2 I\\
			\vdots& & \ddots &  \vdots& \vdots\\
			0  &0 & \cdots&A[\lambda_n+t]& \beta_n I\\
			-\beta_1 I & -\beta_2 I & \cdots& -\beta_n I& I
		\end{array}
		\right).
	\end{align*}
	Let
	\begin{align}
		{\cal R}^{(1)}=\Big\{t \ : \ t\in {\cal R} \ \mbox{and} \  A-(t+\lambda_i)I \ \mbox{is nonsingular for $1\leq i \leq n$}\Big\}.
		\label{1212-1}
	\end{align}
	So for $t\in {\cal R}^{(1)}$,  $A[\lambda_i+t]$ is nonsingular for $1\leq i\leq n$. By the Schur complement, for $t\in {\cal R}^{(1)}$,
	\begin{align*}
			  \det\big(C(t)\big)&=\det\Big(A[\lambda_1+t]\Big) \times\det\Big( A[\lambda_2+t] \Big) \times\cdots \times  \det\Big( A[\lambda_n+t] \Big)\\
			& \quad  \times  \det\Big(I+\beta_1^2 (A[\lambda_1+t])^{-1}+\beta_2^2  (A[\lambda_2+t])^{-1}+\cdots +\beta_n^2 (A[\lambda_n+t])^{-1} \Big)\\
			&= \det\Big(A^{\top}-(t+\lambda_n)I \Big)\times \cdots\times \det\Big(A^{\top}-(t+\lambda_2)I \Big)\times \det\Big(A^{\top}-(t+\lambda_1)I \Big)\\
			& \quad \times \det\Big(I+\beta_1^2(A^{\top}-(t+\lambda_1) I)^{-1}\cdot(A-(t+\lambda_1) I)^{-1}\\
			&\quad\quad\quad+\cdots + \beta_n^2(A^{\top}-(t+\lambda_n) I)^{-1}\cdot (A-(t+\lambda_n) I)^{-1} \Big)\\
			& \quad \times \det\Big(A-(t+\lambda_1)I \Big)\times \det\Big(A-(t+\lambda_2)I \Big)\times\cdots\times \det\Big(A-(t+\lambda_n)I)\Big) \\
			&= \det\Big((A^{\top}-(t+\lambda_n)I)\times \cdots\times(A^{\top}-(t+\lambda_2)I)\cdot(A^{\top}-(t+\lambda_1)I) \Big)\\
			& \quad \times\det\Big(I+\beta_1^2(A^{\top}-(t+\lambda_1) I)^{-1}\cdot(A-(t+\lambda_1) I)^{-1}\\
			&\quad\quad\quad+\cdots + \beta_n^2(A^{\top}-(t+\lambda_n) I)^{-1}(A-(t+\lambda_n) I)^{-1} \Big)\\
			& \quad \times \det\Big((A-(t+\lambda_1)I)\cdot(A-(t+\lambda_2)I)\times\cdots\times(A-(t+\lambda_n)I) \Big).
		\end{align*} 
	Note that for $1\leq i,j\leq n$,
	\begin{align*}
		&\Big(A-(t+\lambda_i)I\Big)\Big(A-(t+\lambda_j)I\Big)=\Big(A-(t+\lambda_j)I\Big)\Big(A-(t+\lambda_i)I\Big);\\
		&\Big(A-(t+\lambda_i)I\Big)^{-1}\Big(A-(t+\lambda_j)I\Big)=\Big(A-(t+\lambda_j)I\Big)\Big(A-(t+\lambda_i)I\Big)^{-1}.
	\end{align*}
	Then, for $t\in {\cal R}^{(1)}$,
	\begin{align*}
		 \det\big(C(t)\big)=& \det\Big( \Big(\Pi_{\ell=1}^n (A^\top -(t+\lambda_\ell)I)\Big) \cdot \Big(\Pi_{\ell=1}^n (A-(t+\lambda_\ell)I)\Big)\\
		& +\sum_{i=1}^n \beta^2_i \Big(\Pi_{\ell\neq i,\ell =1}^n (A^\top -(t+\lambda_\ell)I)\Big) \cdot \Big(\Pi_{\ell\neq i, \ell=1}^n (A-(t+\lambda_\ell)I)\Big) \Big).
	\end{align*}
Recall that $\beta_i \neq 0$ and $A-(t+\lambda_i)I$ is nonsingular. Then each term in above equation is positive definite. By the Minkowski inequality, for $t\in {\cal R}^{(1)}$,
\begin{align}
	 \det\big(C(t)\big)&\geq   \det\Big( \sum_{i=1}^n \beta^2_i \Big(\Pi_{\ell\neq i,\ell =1}^n (A^\top -(t+\lambda_\ell)I)\Big) \cdot \Big(\Pi_{\ell\neq i, \ell=1}^n (A-(t+\lambda_\ell)I)\Big)
	 \Big) \nonumber \\
	 &\triangleq  \det\big(D(t)\big).\label{1212-10}
\end{align}
Let
\begin{align*}
	A_{i,k}(t)&=\Big(A-(t+\lambda_1) I\Big)\cdots\Big(A-(t+\lambda_n) I\Big)\\
	& \ \ \ \cdot \Big(\Big(A-(t+\lambda_{i}) I\Big) \cdots\Big(A-(t+\lambda_{[i+k-1]_n}) I\Big)\Big)^{-1},
\end{align*}
for $1\leq i, k \leq n$. That is, $A_{i,k}(t)$ is the product of the whole terms $A-(\lambda_\ell+t)I$ with $\ell\notin \{i,\ldots,[i+k-1]_n\}$.   Then
\begin{align*}
	 \det\big(D(t)\big)= \det\Big(\beta_1^2 A_{1,1}^{\top}(t)A_{1,1}(t)+\beta_2^2 A_{2,1}^{\top}(t)A_{2,1}(t)+\cdots+ \beta_n^2 A_{n,1}^{\top}(t)A_{n,1}(t) \Big).
\end{align*}
Equivalently,
\begin{align}
	2^n \cdot \det\big(D(t)\big)= \det\Big(2\beta_1^2 A_{1,1}^{\top}(t)A_{1,1}(t)+2\beta_2^2 A_{2,1}^{\top}(t)A_{2,1}(t)+
	\cdots +2\beta_n^2 A_{n,1}^{\top}(t)A_{n,1}(t) \Big).\label{1210-2}
\end{align}
Now we put the two adjacent terms together and take out the common factor, that is, for $1\leq i \leq n$,
\begin{align}
	&\beta_i^2 A_{i,1}^{\top}(t)A_{i,1}(t)+\beta_{[i+1]_n}^2 A_{[i+1]_n,1}^{\top}(t)A_{[i+1]_n,1}(t)\nonumber\\
	&=A_{i,2}^{\top}(t)\left[\beta_i^2\Big(A^{\top}-(t+\lambda_{[i+1]_n}) I\Big)\Big(A-(t+\lambda_{[i+1]_n} )I\Big)\right.\nonumber\\
	& \ \ \ \ \ \ \ \ \ \ \ \ \  \left. +\beta_{[i+1]_n}^2\Big (A^{\top}-(t+\lambda_i) I\Big)
	\Big(A-(t+\lambda_i) I\Big) \right]A_{i,2}(t).\label{1212-2}
\end{align}
Thus, after $\beta_i^2 A_{i,1}^{\top}(t)A_{i,1}(t)$ is coupled with $\beta_{i+1}^2 A_{i+1,1}^{\top}(t)A_{i+1,1}(t)$ for $1\leq i\leq n-1$,
and  $\beta_n^2 A_{n,1}^{\top}(t)A_{n,1}(t)$ is coupled with $\beta_{1}^2 A_{1,1}^{\top}(t)A_{1,1}(t)$,  by (\ref{1210-2}),
\begin{align*}
	2^n\cdot \det \big(D(t)\big)&= \det\Big( A_{1,2}^{\top}(t)\Big[\beta_1^2\Big(A^{\top}-(t+\lambda_{2}) I\Big)\Big(A-(t+\lambda_{2}) I\Big)\\
	& \ \ \ \ \ \ \ \ \ \ \ \ \ +\beta_{2}^2
	\Big(A^{\top}-(t+\lambda_1) I\Big)\Big(A-(t+\lambda_1) I\Big) \Big]A_{1,2}(t) \\
	& \ \ \ +A_{2,2}^{\top}(t)\Big[\beta_2^2\Big(A^{\top}-(t+\lambda_{3}) I\Big)\Big(A-(t+\lambda_{3}) I\Big)\\
	& \ \ \ \ \ \ \ \ \ \ \ \ \ \ \ +\beta_{3}^2\Big (A^{\top}-(t+\lambda_2) I\Big)\Big(A-(t+\lambda_2) I\Big) \Big]A_{2,2}(t)\\
	& \ \ \ +\cdots\\
	& \ \ \ +A_{n-1,2}^{\top}(t)\Big[\beta_{n-1}^2\Big(A^{\top}-(t+ \lambda_{n}) I\Big)\Big(A-(t+ \lambda_{n}) I\Big)\\
	& \ \ \ \ \ \ \ \ \ \ \ \ \ \ \ \ \ +\beta_{n}^2\Big (A^{\top}-(t+ \lambda_{n-1}) I\Big)\Big(A-(t+ \lambda_{n-1}) I\Big) \Big]A_{n-1,2}(t)\\
	& \ \ \ +A_{n,2}^{\top}(t)\Big[\beta_n^2\Big(A^{\top}-(t+\lambda_{1}) I\Big)\Big(A-(t+\lambda_{1}) I\Big)\\
	& \ \ \ \ \ \ \ \ \ \ \ \ \ \ \  +\beta_{1}^2 \Big(A^{\top}-(t+\lambda_n) I\Big)\Big(A-(t+\lambda_n) I\Big) \Big]A_{n,2}(t)
	 \Big).
\end{align*}
Observe that
\begin{align}
	&\beta_i^2\Big(A^{\top}-(t+\lambda_{j}) I\Big)\Big(A-(t+\lambda_{j}) I\Big)+\beta_{j}^2 \Big(A^{\top}-(t+\lambda_i) I\Big)
	\Big(A-(t+\lambda_i )I\Big)\nonumber\\
	& \ \ \ =(\beta_i^2+\beta_j^2)A^{\top}A-\Big[\beta_i^2(t+\lambda_j)+\beta_j^2(t+\lambda_i)\Big](A^{\top}+A)+\Big[\beta_i^2(t+\lambda_j)^2+\beta_j^2(t+\lambda_i)^2\Big]I\nonumber\\
	& \ \ \ =\bigg(\sqrt{\beta_i^2+\beta_j^2}A^{\top}-\frac{\beta_i^2(t+\lambda_j)+\beta_j^2(t+\lambda_i)}{\sqrt{\beta_i^2+\beta_j^2}}I \bigg)\nonumber\\
	& \ \ \ \ \ \ \ \ \times \bigg(\sqrt{\beta_i^2+\beta_j^2}A-\frac{\beta_i^2(t+\lambda_j)+\beta_j^2(t+\lambda_i)}{\sqrt{\beta_i^2+\beta_j^2}}I \bigg)  +\frac{\beta_i^2\beta_j^2(\lambda_i-\lambda_j)^2}{\beta_i^2+\beta_j^2}I. \label{1212-3}
\end{align}
Let
\begin{align*}
	{\cal R}^{(2)}&=\Big\{ t \ : \ t\in {\cal R}\ \mbox{and} \
	\sqrt{\beta_i^2+\beta_{[i+1]_n}^2}A-\frac{\beta_i^2(t+\lambda_{[i+1]_n})+\beta_{[i+1]_n}^2(t+\lambda_i)}{\sqrt{\beta_i^2+\beta_{[i+1]_n}^2}}I \\
	& \ \ \ \ \ \ \ \ \ \ \ \ \mbox{is  nonsingular for} \ 1\leq i \leq n\Big\}.
\end{align*}
By the Minkowski inequality again, for $t\in {\cal R}^{(1)}\cap {\cal R}^{(2)}$,
\begin{align*}
	2^n  \det\big(D(t)\big) \geq  \det\Big( \gamma_{1}^{(2)}A_{1,2}^{\top}(t)A_{1,2}(t)	+\gamma_{2}^{(2)}A_{2,2}^{\top}(t) A_{2,2}(t)+\cdots+\gamma_{n}^{(2)}A_{n,2}^{\top}(t)A_{n,2}(t)
	 \Big).
\end{align*}
Equivalently,
\begin{align}
	2^{2n}  \det\big(D(t)\big) \geq  \det\Big( 2\gamma_{1}^{(2)}A_{1,2}^{\top}(t)A_{1,2}(t)	+2\gamma_{2}^{(2)}A_{2,2}^{\top}(t) A_{2,2}(t)+\cdots+2\gamma_{n}^{(2)}A_{n,2}^{\top}(t)A_{n,2}(t)
	 \Big).\label{1212-4}
\end{align}
Similar to (\ref{1212-2}), again put the two adjacent terms together and take out the common factor, that is, for $1\leq i \leq n$,
\begin{align}
	&\gamma_i^{(2)}A_{i,2}^{\top}(t)A_{i,2}(t)+\gamma_{[i+1]_n}^{(2)}A_{[i+1]_n,2}^{\top}(t)A_{[i+1]_n,2}(t)\nonumber\\
	& \ \ \ =A_{i,3}^{\top}(t)\bigg[\gamma_{i}^{(2)} \Big(A^{\top}-(t+\lambda_{[i+2]_n}) I\Big)\Big(A-(t+\lambda_{[i+2]_n}) I\Big)\nonumber\\
	& \ \ \ \ \ \ \ \ \ \ \ \ \ \ \ +
	\gamma_{[i+1]_n}^{(2)}\Big(A^{\top}-(t+\lambda_{i}) I\Big)\Big(A-(t+\lambda_{i}) I\Big)\bigg]A_{i,3}(t).\label{1212-5}
\end{align}
Similar to (\ref{1212-3}),
\begin{align}
	&\gamma_{i}^{(2)} \Big(A^{\top}-(t+\lambda_{[i+2]_n}) I\Big)\Big(A-(t+\lambda_{[i+2]_n} )I\Big)+
	\gamma_{[i+1]_n}^{(2)}\Big(A^{\top}-(t+\lambda_{i}) I\Big)\Big(A-(t+\lambda_{i})I\Big)\nonumber\\
	& \ \ \  = \bigg(\sqrt{\gamma_{i}^{(2)}+\gamma_{[i+1]_n}^{(2)}}A^{\top}-\frac{\gamma_{i}^{(2)}(t+\lambda_{[i+2]_n})+
		\gamma_{[i+1]_n}^{(2)}(t+\lambda_i)}{\sqrt{\gamma_{i}^{(2)}+\gamma_{[i+1]_n}^{(2)}}}I \bigg)\label{1212-6}\\
	& \ \ \ \quad \times \bigg(\sqrt{\gamma_{i}^{(2)}+\gamma_{[i+1]_n}^{(2)}}A-\frac{\gamma_{i}^{(2)}(t+\lambda_{[i+2]_n})+\gamma_{[i+1]_n}^{(2)}(t+\lambda_i)}
	{\sqrt{\gamma_{i}^{(2)}+\gamma_{[i+1]_n}^{(2)}}}I \bigg) +\frac{\gamma_{i}^{(2)}\gamma_{[i+1]_n}^{(2)}(\lambda_i-\lambda_{[i+2]_n})^2}{\gamma_{i}^{(2)}+\gamma_{[i+1]_n}^{(2)}}I.
	\nonumber
\end{align}
Let
\begin{align*}
	{\cal R}^{(3)}=\Big\{ t  & : \ t\in {\cal R}\ \mbox{and} \  \sqrt{\gamma_{i}^{(2)}+\gamma_{[i+1]_n}^{(2)}}A-\frac{\gamma_{i}^{(2)}(t+\lambda_{[i+2]_n})+\gamma_{[i+1]_n}^{(2)}(t+\lambda_i)}{\sqrt{\gamma_{i}^{(2)}
			+\gamma_{[i+1]_n}^{(2)}}}I \\
	& \ \ \ \ \ \mbox{is nonsingular for $1\leq i\leq n$}\Big\}.
\end{align*}
We have, by the Minkowski inequality, that for $t\in \cap_{\ell=1}^3{\cal R}^{(\ell)}$,
\begin{align}
	2^{2n}  \det\big(D(t)\big)  \geq  \det\Big( \sum_{i=1}^n \gamma_{i}^{(3)} A_{i,3}^{\top}(t)A_{i,3}(t) \Big).\label{1212-7}
\end{align}
In general, for $k=3,\ldots, n-2$, define
\begin{align*}
	{\cal R}^{(k+1)}=\Big\{ t  & : \ t\in {\cal R}\ \mbox{and} \  \sqrt{\gamma_{i}^{(k)}+\gamma_{[i+1]_n}^{(k)}}A-\frac{\gamma_{i}^{(k)}(t+\lambda_{ [i+k]_n})+\gamma_{[i+1]_n}^{(k)}(t+\lambda_i)}{\sqrt{\gamma_{i}^{(k)}
			+\gamma_{[i+1]_n}^{(k)}}}I \\
	& \ \ \ \ \ \mbox{is nonsingular for $1\leq i\leq n$}\Big\},
\end{align*}
and
\begin{align*}
	{\cal R}^{(n)}=\Big\{ t  & : \ t\in {\cal R}\ \mbox{and} \
	\sqrt{\sum_{i=1}^n \gamma_{i}^{(n-1)}}A^{\top}-\frac{\sum_{i=1}^n \gamma_{i}^{(n-1)}(t+\lambda_{[n+i-1]_n})}{\sqrt{\sum_{i=1}^n \gamma_{i}^{(n-1)}}}I
	\\
	& \ \ \ \ \ \mbox{is nonsingular for $1\leq i\leq n$}\Big\}.
\end{align*}
Repeating  the above procedure given by (\ref{1212-4})-(\ref{1212-7}), we get that for $t\in \cap_{\ell=1}^{n-1}{\cal R}^{(\ell)}$,
\begin{align*}
	2^{(n-2)n}  \det\big(D(t)\big)
	&\geq  \det\bigg( \sum_{i=1}^n \gamma_{i}^{(n-1)} A_{i,n-1}^{\top}(t)A_{i,n-1}(t)             \bigg)\\
	&= \det\bigg( \sum_{i=1}^n \gamma_{i}^{(n-1)} \Big(A^{\top}-(t+\lambda_{[n+i-1]_n})I\Big)\Big(A-(t+\lambda_{[n+i-1]_n})I\Big)            \bigg)\\
	&= \det\bigg( \sum_{i=1}^n \gamma_{i}^{(n-1)}A^{\top}A- \sum_{i=1}^n \gamma_{i}^{(n-1)}(t+\lambda_{[n+i-1]_n})(A^{\top}+A) \\
	& \quad \ \ + \sum_{i=1}^n \gamma_{i}^{(n-1)}(t+\lambda_{[n+i-1]_n})^2I
	 \bigg)\\
	&= \det\bigg(\bigg( \sqrt{\sum_{i=1}^n \gamma_{i}^{(n-1)}}A^{\top}-\frac{\sum_{i=1}^n \gamma_{i}^{(n-1)}(t+\lambda_{[n+i-1]_n})}{\sqrt{\sum_{i=1}^n \gamma_{i}^{(n-1)}}}I
	\bigg) \\
	&  \quad  \ \ \times \bigg( \sqrt{\sum_{i=1}^n \gamma_{i}^{(n-1)}}A-\frac{\sum_{i=1}^n \gamma_{i}^{(n-1)}(t+\lambda_{[n+i-1]_n})}
	{\sqrt{\sum_{i=1}^n \gamma_{i}^{(n-1)}}}I
	\bigg)\\
	& \quad +\bigg(\sum_{i=1}^n \gamma_{i}^{(n-1)}(t+\lambda_{[n+i-1]_n})^2-\frac{\Big(\sum_{i=1}^n \gamma_{i}^{(n-1)}(t+\lambda_{[n+i-1]_n})\Big)^2}{\sum_{i=1}^n \gamma_{i}^{(n-1)}}
	\bigg)I \bigg).
\end{align*}
It can be directly calculated that
\begin{align*}
	&\sum_{i=1}^n \gamma_{i}^{(n-1)}(t+\lambda_{[n+i-1]_n})^2-\frac{\Big(\sum_{i=1}^n \gamma_{i}^{(n-1)}(t+\lambda_{[n+i-1]_n})\Big)^2}{\sum_{i=1}^n \gamma_{i}^{(n-1)}}\\
	& \ \ \ =\frac{\sum_{i\neq j}\gamma_{i}^{(n-1)}\gamma_{j}^{(n-1)}(\lambda_{[n+i-1]_n}-\lambda_{[n+j-1]_n})^2}{\sum_{i=1}^n \gamma_{i}^{(n-1)}}.
\end{align*}
By  the Minkowski inequality, we have that for $t\in \cap_{\ell=1}^n{\cal R}^{(\ell)}$,
\[
 \det\big(D(t)\big) \geq  2^{-n(n-2)}\bigg(\frac{\sum_{i\neq j}\gamma_{i}^{(n-1)}\gamma_{j}^{(n-1)}(\lambda_{[n+i-1]_n}-\lambda_{[n+j-1]_n})^2}
{\sum_{i=1}^n \gamma_{i}^{(n-1)}}\bigg)^n.
\]
As the set $\cap_{\ell=1}^n{\cal R}^{(\ell)}$ is dense in ${\cal R}$, by (\ref{1212-9}), (\ref{1212-10}) and the continuity of $ \det\big(C(t)\big)$ in $t$, we have
\begin{align*}
	 \det(BB^\top)=\lim_{t\rightarrow 0}  \det\big(C(t)\big) &\geq \lim_{t\rightarrow 0} \det\big(D(t)\big) \\
	& \geq  2^{-n(n-2)}\bigg(\frac{\sum_{i\neq j}\gamma_{i}^{(n-1)}\gamma_{j}^{(n-1)}(\lambda_{[n+i-1]_n}-\lambda_{[n+j-1]_n})^2}
	{\sum_{i=1}^n \gamma_{i}^{(n-1)}}\bigg)^n.
\end{align*}
This completes the proof of the  proposition.
\hfill$\Box$\\[-0.05in]

With Proposition \ref{tech-prop} in hand,  we establish (\ref{uniform_bound}).
\begin{myprop}\label{p-bound}
	If $P\in {\cal R}^{n\times n}$ satisfies {\rm (\ref{add_con})}, then $\|P\|\leq {\rho}$ with
	\begin{align*}
		{\rho}&=2^{n(n-2)}\bigg(\frac{\sum_{i=1}^n \gamma_{i}^{(n-1)}}{\sum_{i\neq j}\gamma_{i}^{(n-1)}\gamma_{j}^{(n-1)}(\lambda_{[n+i-1]_n}-\lambda_{[n+j-1]_n})^2}\bigg)^n\Big (\|\beta\|^2+4\Big)^{(n^2-1)}\\
		& \ \ \ \times (n^2-1)! (n+\lambda_1+\|\beta\|) \max\{\sqrt{\bar {\rho}},1\},
	\end{align*}
	where $\bar {\rho}$ is defined in {\rm (\ref{1212-37})}, and $\gamma^{(n-1)}_i$ $(1\leq i\leq n)$ are given by {\rm (\ref{1212-52})}.
\end{myprop}

\noindent
{\it Proof}: \  Consider any $Y\in {\cal R}^{(n+1)\times n}$ with
\begin{align}
	\left(\begin{array}{ccccc}
		-\sqrt{\bar {\rho}}{\bf 1}^\top\\
		-\sqrt{\bar {\rho}}{\bf 1}^\top\\
		\vdots\\
		-\sqrt{\bar {\rho}}{\bf 1}^\top\\
		{\bf 0}^\top
	\end{array}\right)\leq
	\left(\begin{array}{ccccc}
		Y_{(1)}^{\top}\\
		Y_{(2)}^{\top}\\
		\vdots\\
		Y_{(n)}^{\top}\\
		Y_{(n+1)}^{\top}
	\end{array}\right)
	\leq \left(\begin{array}{ccccc}
		\sqrt{\bar {\rho}}{\bf 1}^\top\\
		\sqrt{\bar {\rho}}{\bf 1}^\top\\
		\vdots\\
		\sqrt{\bar {\rho}}{\bf 1}^\top\\
		{\bf 1}^\top
	\end{array}\right)
	\ \mbox{and} \ B^\top \left(\begin{array}{cccc}
		P_{(1)}^{\top}\\
		P_{(2)}^{\top}\\
		\vdots\\
		P_{(n)}^{\top}
	\end{array}\right)=\left(\begin{array}{cccc}
		Y_{(1)}^{\top}\\
		Y_{(2)}^{\top}\\
		\vdots\\
		Y_{(n+1)}^{\top}
	\end{array}\right).
	\label{1212-33}
\end{align}
It suffices to show that $P$ in the last equation of (\ref{1212-33}) satisfies $\|P\|\leq {\rho}$.
To this end, we first have
\[
B B^\top \left(\begin{array}{cccc}
	P_{(1)}^{\top}\\
	P_{(2)}^{\top}\\
	\vdots\\
	P_{(n)}^{\top}
\end{array}\right)=B\left(\begin{array}{cccc}
	Y_{(1)}^{\top}\\
	Y_{(2)}^{\top}\\
	\vdots\\
	Y_{(n+1)}^{\top}
\end{array}\right).
\]
Hence,
\begin{align}
	\left(\begin{array}{cccc}
		P_{(1)}^{\top}\\
		P_{(2)}^{\top}\\
		\vdots\\
		P_{(n)}^{\top}
	\end{array}\right)=\Big( B B^\top\Big)^{-1} B\left(\begin{array}{cccc}
		Y_{(1)}^{\top}\\
		Y_{(2)}^{\top}\\
		\vdots\\
		Y_{(n+1)}^{\top}
	\end{array}\right)=\frac{1}{ \det(BB^\top)} \textrm{adj} ( B B^\top) B\left(\begin{array}{cccc}
		Y_{(1)}^{\top}\\
		Y_{(2)}^{\top}\\
		\vdots\\
		Y_{(n+1)}^{\top}
	\end{array}\right).
	\label{1212-34}
\end{align}
Using (\ref{1212-11}),
\begin{align*}
	\|(A-\lambda_iI)(A^{\top}-\lambda_iI)\|=\|AA^{\top}-\lambda_i(A+A^{\top})+\lambda_i^2I\|\leq 4.
\end{align*}
According to the definition of $B^\top$ given by (\ref{12-10-1}), it follows from $C(0)=BB^\top$ and (\ref{1212-50}) that
\begin{align*}
	\|BB^{\top}\|\leq \|\beta\|^2+4.
\end{align*}
Since each entry of $ \textrm{adj}(BB^{\top})$ is  an $(n^2-1)\times (n^2-1)$ determinant and it has $(n^2-1)!$ terms,
\begin{align}
	\|\textrm{adj}(BB^{\top}) \| \leq \Big (\|\beta\|^2+4\Big)^{(n^2-1)}(n^2-1)!.
	\label{1212-35}
\end{align}
By (\ref{1212-33}) and the definition of $B^\top$ given by (\ref{12-10-1}),
\begin{align}
	\left\| B\left(\begin{array}{cccc}
		Y_{(1)}^{\top}\\
		Y_{(2)}^{\top}\\
		\vdots\\
		Y_{(n+1)}^{\top}
	\end{array}\right)\right\|
	\leq (n+\lambda_1+\|\beta\|) \max\{\sqrt{\bar {\rho}},1\}.
	\label{1212-36}
\end{align}
Using Proposition \ref{tech-prop}, (\ref{1212-34})-(\ref{1212-36}), we know the proposition holds.\hfill$\Box$\\[-0.05in]

As any optimal solution to {\bf (OP[A])} given by (\ref{july-15-4}) satisfies \eqref{add_con} regardless of how to choose $A \in {\cal R}^{n\times n}_+$ with $A{\bf 1}^\top\leq {\bf 1}^\top$, using Propositions 1 and 2, we have (\ref{uniform_bound}). This is formally summarized by the following proposition.
\begin{myprop}\label{tech-prop-bound}
	For any $A \in {\cal R}^{n\times n}_+$ with $A{\bf 1}^\top\leq {\bf 1}^\top$, any optimal solution to {\bf (OP[A])} given by {\rm (\ref{july-15-4})}, denoted by $P_{\sf sln}(A)$, must satisfy $\|P_{\sf sln}(A)\|\leq {\rho}$  with ${\rho}$ being given by Proposition {\rm \ref{p-bound}}.
\end{myprop}

Now we can get the effectiveness of the {\bf AM} for the nonconvex optimization problem {\bf (OP)} given by (\ref{july-15-2}).
\begin{mythm}\label{KKT}
	The algorithm $($given by {\bf Step-a} and {\bf Step-b}$)$ can guarantee$:$
	\begin{itemize}
		\item[{\rm (i)}]  The sequence $\{{\cal L}(A_\ell,P_\ell):\ell\geq 0\}$  has a finite limit denoted by ${\cal L}_\infty;$
		\item[{\rm (ii)}] The sequence $\{(A_\ell, P_\ell): \ell\geq 0\}$ admits at least one limit point, and its any limit point, say $( {A}_\infty,{P}_\infty)$, must  satisfy {\rm (\ref{1212-40})-(\ref{1212-31})}, the KKT  conditions of the nonconvex optimization problem $(\bf OP)$ given by {\rm (\ref{july-15-2})}. Further, $${\cal L}( {A}_\infty,{P}_\infty)= {\cal L}_\infty;$$ \setstretch{0.01}
		\item[{\rm (iii)}]  Any limit point admitted by $\{(A_\ell, P_\ell): \ell\geq 0\}$ must not be local maximum.
	\end{itemize}
\end{mythm}

\begin{myrmk}\label{rem-critical-point}
	{\rm Under the assumption in which the sequence $\{(A_\ell, P_\ell): \ell\geq 0\}$ admits a limit point,
		with the help of the asymptotically critical point result given by {\rm \cite{bertsekas1995}}, {\rm\cite{grippo2000convergence}} proves the admitted limit point must be a critical point by the contradiction method. Recall that the critical points and the solutions to the KKT conditions are  the same, see Theorem 12.1 in \cite{Nocedal-2006}.  In addition to establishing the KKT  conditions hold for the admitted limit point, Theorem {\rm \ref{KKT} (ii)} not only confirms that the existence of the limit point of $\{(A_\ell, P_\ell): \ell\geq 0\}$ but also claims that all the limit points give the same objective value ${\cal L}_\infty$. Different from the method employed by
		{\rm\cite{grippo2000convergence}}, here we directly prove any limit point admitted by the sequence $\{(A_\ell, P_\ell): \ell\geq 0\}$ satisfies the KKT  conditions by establishing the boundedness of the Lagrange multipliers.}
\end{myrmk}
	\begin{myrmk}\label{rem-sequence-convergence}
		{\rm Although Theorem {\rm \ref{KKT} (ii)} gives that the existence of the limit point of $\{(A_\ell, P_\ell): \ell\geq 0\}$, the sequence $\{(A_\ell, P_\ell): \ell\geq 0\}$ itself cannot  be proved to be convergent. The main difficulty to prove its convergence is that we cannot exclude the possibility that there exists a solution to $(\bf{OP[A]})$ at some iteration step such that $P_{\hat\ell}$ is singular. As long as $P_{\hat \ell}$ is singular, the solution to $(\bf{OP[P]})$ with ${\bf P}=P_{\hat \ell}$ may not be unique because (i) there exists a matrix $Q \in {\cal R}^{n\times n}$ such that $P_{\hat \ell}Q={\bf 0}$;
			(ii) for such $Q$, ${\cal L}(A, P_{\hat \ell})$ and ${\cal L}(A+tQ, P_{\hat \ell})$ are equal with all $t\in {\cal R}_+$; and (iii) possibly for small enough $t$, $A+tQ$ is nonnegative and $(A+tQ){\bf 1}^\top\leq {\bf 1}^\top$
			when $A$ is a solution to $(\bf{OP[P]})$ with ${\bf P}=P_{\hat \ell}$. Thus, the singularity of $P_{\hat\ell}$ incurs the nonuniqueness for the solution $A_{\hat\ell+1}$ to $(\bf{OP[P_{\hat\ell}]})$. The nonuniqueness feature of $A_{\hat\ell+1}$ breaks down the existence of the limit of the entire sequence $\{(A_\ell, P_\ell): \ell\geq 0\}$, see Example \ref{examp-noconverge} for a concrete case.
			At the same time,  Theorem \ref{KKT} can only ensure the algorithm finds the solutions to the KKT conditions instead of the global optimal solutions. }
	\end{myrmk}

\noindent
{\it Proof} of {\bf Theorem \ref{KKT}}: We first prove that $\{{\cal L}(A_\ell, P_\ell): \ell\geq 0\}$ is monotonically decreasing.
Recalling the procedure to generate $\{(A_\ell, P_\ell): \ell\geq 0\}$ by {\bf Step-a} and {\bf Step-b}, we first arbitrarily
pick up $A_0$ which satisfies $A_0\in {\cal R}_+^{n\times n}$ and $A_0{\bf 1}^\top\leq {\bf 1}^\top$, and
\begin{align}
	\left\{
	\begin{array}{llll}
		&P_\ell \ \mbox{is a solution for the convex optimization problem $(\bf OP[A])$}\\
		& \ \mbox{given by (\ref{july-15-4}) with $A=A_{\ell}$ for } \ell\geq 0;\\
		& A_\ell \ \mbox{is a solution for the convex optimization problem $(\bf{OP[P]})$}\\
		& \ \mbox{given by (\ref{july-15-6}) with $P=P_{\ell-1}$ for }   \ell\geq 1.
	\end{array}
	\right.\label{1212-59}
\end{align}	
Hence we have that for $\ell\geq 0$,
\begin{align}
	{\cal L}(A_\ell, P_\ell)\geq {\cal L}(A_{\ell+1}, P_\ell) \geq {\cal L}(A_{\ell+1}, P_{\ell+1}). \label{limit-4}
\end{align}	
Since ${\cal L}(A_\ell, P_\ell) \geq 0$ for each $\ell$, there exists a nonnegative real number $ {\cal L}_\infty$ such that
\begin{eqnarray}
	\lim_{\ell \rightarrow \infty} {\cal L}(A_\ell, P_\ell) ={\cal L}_\infty.\label{limit-1}
\end{eqnarray}
This gives the first part of the theorem.

Now prove the second part of the theorem. The existence of the limit point admitted by the sequence $\{(A_\ell, P_\ell): \ell\geq 0\}$ directly follows from
Proposition \ref{tech-prop-bound}.
Further, by the first part of the theorem, we know that the objective value given by any limit point of
$\{(A_\ell, P_\ell): \ell\geq 0\}$ must be ${\cal L}_\infty$. So the remainder of the proof of the second part is to show that
any limit point of  $\{(A_\ell, P_\ell): \ell\geq 0\}$  must satisfy (\ref{1212-40})-(\ref{1212-31}).

To this end, for any convergent subsequence of $\{(A_{\ell_k},P_{\ell_k}): k\geq 1\}$, let $(A_\infty,P_\infty)$ be the corresponding limit.
It suffices to show that
\begin{align*}
		(A_\infty,P_\infty) \ \ \mbox{satisfies \ \ (\ref{1212-40})-(\ref{1212-31})}.
\end{align*}
$P_{\ell_k}$ is a solution to the convex optimization problem ({\bf OP[A]}) given by (\ref{july-15-4}) with $A=A_{\ell_k}$.
Hence, $(A_{\ell_k},P_{\ell_k})$ satisfies the {\bf KKT} condition given by (\ref{1212-40}) for each $k\geq 1$.
That is, for each $k\geq 1$, there exist
$\{\alpha_{\ell_k,i}: 1\leq i\leq n\}$ and $\{\widetilde \alpha_{\ell_k,i}: 1\leq i\leq n\}$ such that
\begin{eqnarray}
	\left\{
	\begin{array}{lllllll}
		&\frac{\partial}{\partial p_{ij}} {\cal L}(A_{\ell_k},P_{\ell_k})+\alpha_{\ell_k,i} -\widetilde\alpha_{\ell_k,j}\cdot \beta_i=0, \ 1\leq i,j\leq n;\\
		& P_{\ell_k}{\bf 1}^\top={\bf 1}^\top;\\
		& \beta P_{\ell_k} \geq {\bf 0};\\
		& \widetilde\alpha_{\ell_k,j} \cdot \Big(\sum_{i=1}^n \beta_i p_{\ell_k,ij}\Big)=0, \ 1\leq j\leq n;\\
		& (\widetilde\alpha_{\ell_k,1},\ldots, \widetilde\alpha_{\ell_k,n})\in {\cal R}^n_+.
	\end{array}
	\right.
	\label{1212-53}
\end{eqnarray}
To show $(A_\infty,P_\infty)$ satisfies (\ref{1212-40})-(\ref{1212-31}), we work on the sequence  $(A_{\ell_k},P_{\ell_k})$.
In view of the equations except the second and third ones in (\ref{1212-53}), we first establish the boundedness of two sequences
	$\{(\alpha_{\ell_k,i}, 1\leq i\leq n): k\geq 1\}$ and $\{(\widetilde \alpha_{\ell_k,i}, 1\leq i\leq n): k\geq 1\}$.
Let
\[
\widetilde {\cal O}_{\ell_k}=\Big\{j: \ \sum_{i=1}^n \beta_i {p}_{\ell_k,ij}>0\Big\}.
\]
From the boundedness of $\{(A_{\ell_k},P_{\ell_k}): k\geq 1\}$ guaranteed by Proposition \ref{tech-prop-bound}, as $n$ is finite, there exist a subset of $\{1,\ldots,n\}$ denoted by $\widetilde {\cal O}$ and a subsequence of $\{\ell_k: k\geq 1\}$, denoted by $\{\ell_k: k\geq 1\}$ for the sake of the notation simplicity,  such that for $k\geq 1$,
\begin{align}
	\widetilde {\cal O}_{\ell_k}=\widetilde {\cal O},  \ \ \lim_{k\rightarrow \infty}A_{\ell_k}={A}_\infty, \ \ \lim_{k\rightarrow \infty}P_{\ell_k}=P_\infty. \label{1212-56}
\end{align}
Then for this subsequence $\{\ell_k: k\geq 1\}$, (\ref{1212-53}) can be written as:
\begin{eqnarray}
	\left\{
	\begin{array}{lllllll}
		&\frac{\partial}{\partial p_{ij}} {\cal L}(A_{\ell_k},P_{\ell_k})+\alpha_{{\ell_k},i} -\widetilde\alpha_{{\ell_k},j}\cdot \beta_i=0, \ 1\leq i,j\leq n;\\
		& P_{\ell_k}{\bf 1}^\top={\bf 1}^\top;\\
		& \beta P_{\ell_k} \geq {\bf 0};\\
		& \sum_{i=1}^n \beta_i p_{{\ell_k},ij}>0 \ \ \mbox{and } \ \widetilde\alpha_{{\ell_k},j}=0, \ j \in \widetilde{\cal O};\\
		& \sum_{i=1}^n \beta_i p_{{\ell_k},ij}=0, \ \ j \notin \widetilde{\cal O};\\
		& (\widetilde\alpha_{\ell_k,1},\ldots, \widetilde\alpha_{\ell_k,n})\in {\cal R}^n_+.
	\end{array}
	\right.
	\label{1212-57}
\end{eqnarray}
Let $\widetilde{n}_o$ be the cardinality of set $\widetilde{\cal O}$.
(\ref{1212-57}) can be considered as an affine mapping from the polyhedral convex set ${\cal R}^n\times{\cal R}_+^{n-\widetilde n_o}$ into ${\cal R}^{n^2}$. Here each component in ${\cal R}^n\times{\cal R}_+^{n-\widetilde n_o}$ corresponds to the Lagrange multipliers $\{\alpha_{\ell_k,i}: 1\leq i\leq n\}$ and $\{(\widetilde \alpha_{\ell_k,i}: i\notin \widetilde{\cal O}\}$. Each image in ${\cal R}^{n^2}$ corresponds to $\frac{\partial}{\partial p_{ij}} {\cal L}(A_{\ell_k},P_{\ell_k})$ $(1\leq i,j\leq n)$.
In view of (\ref{ob_exp}), and the fact that $A_{\ell_k}$ is nonnegative and $A_{\ell_k}{\bf 1}^\top\leq {\bf 1}^\top$, we, by Proposition \ref{tech-prop-bound},  know that there exists a constant $c$ such that
\[
\Big| \frac{\partial}{\partial p_{ij}} {\cal L}(A_{\ell_k},P_{\ell_k})- \frac{\partial}{\partial p_{ij}} {\cal L}(A_{\ell_{\overline k}},P_{\ell_{\overline k}})\Big|
\leq c\cdot \Big(\Big\|P_{\ell_k}-P_{\ell_{\overline k}}\Big\|+\Big\|{A}_{\ell_k}-{A}_{\ell_{\overline k}}\Big\|\Big).
\]
Write
\[
\delta(k,1)=\Big\|P_{\ell_k}-P_{\ell_{1}}\Big\|+\Big\|{A}_{\ell_k}-{A}_{\ell_{1}}\Big\|.
\]
Let $\alpha_{\ell_k}=(\alpha_{\ell_k,1},\ldots,\alpha_{\ell_k,n})$, $\widetilde \alpha_{\ell_k}|_{\widetilde{\cal O}^c}$ be an $(n-\widetilde{n}_o)$-dimensional vector
	given from $\widetilde\alpha_{\ell_k}=(\widetilde\alpha_{\ell_k,1},\ldots,\widetilde\alpha_{\ell_k,n})$ by deleting the coordinates whose indices are in $\widetilde{\cal O}$,
	$J^i\in {\cal R}^{n\times n}$ a matrix whose $i$th column is ${\bf 1}^\top$ while the other entries are zero, and $I_{\widetilde{\cal O}^c}(\beta_i)$ an $n$-by-$(n-\widetilde{n}_o)$ matrix generated from $\beta_iI$ by deleting its columns with indices in $\widetilde{\cal O}$. Then the first, fourth and last equations in (\ref{1212-57}) can be written as
	a linear map ${\cal R}^n\times {\cal R}_+^{n-\widetilde{n}_o}\mapsto {\cal R}^{n^2}$:
	\begin{eqnarray}
		\begin{bmatrix}
			J^1 & -I_{\widetilde{\cal O}^c}(\beta_1)\\
			\vdots & \vdots \\
			J^n & -I_{\widetilde{\cal O}^c}(\beta_n)
		\end{bmatrix}
		\left(
		\begin{array}{cc}
			\alpha_{\ell_k}^\top\\
			\widetilde{\alpha}_{\ell_k}^\top|_{\widetilde{\cal O}^c}
		\end{array}
		\right)
		=-\left(
		\begin{array}{cccc}
			\Big(\frac{\partial}{\partial p_{ij}} {\cal L}(A_{\ell_k},P_{\ell_k})\Big)^\top_{(1)}\\[0.10in]
			\Big(\frac{\partial}{\partial p_{ij}} {\cal L}(A_{\ell_k},P_{\ell_k})\Big)^\top_{(2)}\\
			\vdots\\
			\Big(\frac{\partial}{\partial p_{ij}} {\cal L}(A_{\ell_k},P_{\ell_k})\Big)^\top_{(n)}
		\end{array}
		\right). \label{23-07-1}
	\end{eqnarray}
		As the rank of the coefficient matrix of the left-hand-side of (\ref{23-07-1}) may be less than $(2n-\widetilde{n}_o)$, and ${\cal R}_+^{n-\widetilde{n}_o}$ is not a linear space, we cannot use the generalized inverse (see Section 12.7 in \cite{lachaster}) to get the representation of $\{\alpha_{\ell_k}: k\geq 1\}$ and $\{\widetilde \alpha_{\ell_k}|_{\widetilde{\cal O}^c}: k\geq 1\}$
		in terms of the generalized inverse of the coefficient matrix and $\frac{\partial}{\partial p_{ij}} {\cal L}(A_{\ell_k},P_{\ell_k})$.
		In the following, we use Corollary 7.1 in \cite{lee2005quadratic} to obtain the boundedness of $\{\alpha_{\ell_k}: k\geq 1\}$ and $\{\widetilde \alpha_{\ell_k}|_{\widetilde{\cal O}^c}: k\geq 1\}$. To see Corollary 7.1 in \cite{lee2005quadratic} to be applicable here, first specify its notation, $\Delta$, $A$, $b$, $y$ and $y^\prime$ into what they represent here:
		\begin{eqnarray*}
			&&\Delta= {\cal R}^n\times {\cal R}_+^{n-\widetilde{n}_o}, A=\begin{bmatrix}
				J^1 & -I_{\widetilde{\cal O}^c}(\beta_1)\\
				\vdots & \vdots \\
				J^n & -I_{\widetilde{\cal O}^c}(\beta_n)
			\end{bmatrix}, b={\bf 0}, \\
			&& \mbox{$y$ equals to the right-hand-side of \eqref{23-07-1}},\\
			&& \mbox{$y^\prime$ equals to the right-hand-side of \eqref{23-07-1} with $\ell_k=\ell_1$}.
		\end{eqnarray*}
		This specification shows that the conditions in Corollary 7.1 in \cite{lee2005quadratic} are well satisfied here.
		Then it follows from this corollary that there exist a constant $\overline{c}$ and the Lagrange multipliers $\{\alpha_{\ell_k}: k\geq 2\}$ and $\{\widetilde \alpha_{\ell_k}|_{\widetilde{\cal O}^c}: k\geq 2\}$ (again denoted by this for the sake of the notation simplicity) which satisfy (\ref{1212-57}) and then (\ref{1212-53}) such that for $k\geq 2$,
		\begin{align}
			\Big\|\alpha_{\ell_1} -\alpha_{\ell_k} \Big\|\leq\overline{c}\cdot \delta(k,1) \ \mbox{and} \
			\Big\|\widetilde{\alpha}_{\ell_1}|_{\widetilde{\cal O}^c} - \widetilde{\alpha}_{\ell_k}|_{\widetilde{\cal O}^c} \Big\|\leq\overline{c}\cdot \delta(k,1).\label{1212-62}
		\end{align}
		Hence we find the bounded Lagrange multipliers $\{(\alpha_{\ell_k,i}, 1\leq i\leq n): k\geq 1\}$ and $\{(\widetilde \alpha_{\ell_k,i}, 1\leq i\leq n): k\geq 1\}$ which satisfy (\ref{1212-57}) and then (\ref{1212-53}).
	Therefore, there exists a further subsequence of $\{\ell_k: k\geq 1\}$ denoted by
	$\{\overline{\ell}_k: k\geq 1\}$ such that
	\begin{align*}
		&\lim_{k\rightarrow \infty}  \alpha_{\overline{\ell}_k,i}=\alpha_{\infty,i}, \ \ \ \  \lim_{k\rightarrow \infty}  \widetilde\alpha_{\overline{\ell}_k,i}=\widetilde\alpha_{\infty,i}, \ 1\leq i\leq n.
	\end{align*}
	Further, using (\ref{1212-56}) and taking the limit along the sequence of $\{\overline{\ell}_k:k\geq 1\}$ in (\ref{1212-57}) yields that
	\begin{eqnarray}
		\left\{
		\begin{array}{lllllll}
			&\frac{\partial}{\partial p_{ij}} {\cal L}(A_{\infty},P_{\infty})+\alpha_{\infty,i} -\widetilde\alpha_{{\infty},j}\cdot \beta_i=0, \ 1\leq i,j\leq n;\\
			& P_{\infty}{\bf 1}^\top={\bf 1}^\top;\\
			& \beta P_{\infty} \geq {\bf 0};\\
			& \sum_{i=1}^n \beta_i p_{{\infty},ij}>0 \ \ \mbox{and } \ \widetilde\alpha_{{\infty},j}=0, \ j \in \widetilde{\cal O};\\
			& \sum_{i=1}^n \beta_i p_{{\infty},ij}=0, \ \ j \notin \widetilde{\cal O};\\
			& (\widetilde\alpha_{\infty,1},\ldots, \widetilde\alpha_{\infty,n})\in {\cal R}^n_+.
		\end{array}
		\right.
		\label{1212-58}
	\end{eqnarray}
	If we can show that
	\begin{align}
		\mbox{there exist $\pi_{\infty,i}$ and $\widetilde\pi_{\infty,ij}$ with $1\leq i,j\leq n$ such that (\ref{1212-31}) holds with $(A_\infty,P_\infty)$},
		\label{limit-2}
	\end{align}
	then by (\ref{1212-58}),  $(A_\infty,P_\infty)$ satisfies (\ref{1212-40})-(\ref{1212-31}).
	Thus it remains to show (\ref{limit-2}).
	
	Now let $\overline{A}_{\overline{\ell}_{k+1}}$ be the solution to the convex optimization problem  {\bf (OP[P])} given by (\ref{july-15-6}) with
	$P_0=P_{\overline{\ell}_k}$ for each $k\geq 0$. Since $\{\overline{\ell}_k,k\geq 1\}$ is a subsequence of $\{\ell,\ell\geq 1\}$ and the solutions of {\bf (OP[P])} are not necessarily unique, we know that for each $k\geq 0$, $A_{\overline{\ell}_{k+1}}$ may not be equal to $\overline{A}_{\overline{\ell}_{k+1}}$.
	By (\ref{1212-31}), we have the {\bf KKT} condition for $k\geq 1$,
	\begin{align}	\left\{
		\begin{array}{lllllll}
			&\frac{\partial}{\partial a_{ij}} {\cal L}(\overline{A}_{\overline{\ell}_k},P_{\overline{\ell}_{k-1}})+\pi_{\overline{\ell}_k,i} -\widetilde \pi_{\overline{\ell}_k,ij}=0, \ 1\leq i,j\leq n;\\
			& \overline{A}_{\overline{\ell}_k}{\bf 1}^\top\leq {\bf 1}^\top;\\
			& \pi_{\overline{\ell}_k,i} \times \Big(\sum_{j=1}^n \overline{a}_{\overline{\ell}_k,ij}-1\Big)=0, \ 1\leq i\leq n;\\
			& \overline{A}_{\overline{\ell}_k} \ \mbox{is nonnegative};\\
			& \widetilde \pi_{\overline{\ell}_k,ij}\overline{a}_{\overline{\ell}_k,ij}=0, \ 1\leq i,j\leq n;\\
			&(\pi_{\overline{\ell}_k,1},\ldots,\pi_{\overline{\ell}_k,n})\in {\cal R}^n_+ \ \mbox{and} \ (\widetilde{\pi}_{\overline{\ell}_k,ij})_{n\times n}\in {\cal R}^{n\times n}_+.
		\end{array}
		\right.
		\label{1212-60}
	\end{align}
	Similar to (\ref{1212-62}), we can show the boundedness of two sequences $\{( \pi_{\overline{\ell}_k,i}, 1\leq i\leq n): k\geq 1\}$
	and $\{(\widetilde \pi_{\overline{\ell}_k,ij}, 1\leq i,j\leq n): k\geq 1\}$. Hence, there exists a subsequence of $\{\overline{\ell}_k: k\geq 1\}$, denoted by itself for the sake of the notation simplicity, such that
	\begin{align}
		&\lim_{k\rightarrow \infty}  \pi_{\overline{\ell}_k,i}=\pi_{\infty,i}, \ \ \ \    \lim_{k\rightarrow \infty}  \widetilde\pi_{\overline{\ell}_k,ij}=\widetilde\pi_{\infty,ij}, \ 1\leq i,j\leq n; \ \ \lim_{k\rightarrow \infty} \overline{A}_{\overline{\ell}_k}  =\overline{A}_{\infty}.\label{1212-66}
	\end{align}
	Thus, $\overline{A}_{\infty}$ solves the convex optimization problem  {\bf (OP[P])} given by (\ref{july-15-6}) with
	$P=P_{\infty}$. Furthermore, from (\ref{1212-60}),
	\begin{align}	\left\{
		\begin{array}{lllllll}
			&\frac{\partial}{\partial a_{ij}} {\cal L}(\overline{A}_{\infty}, P_{\infty})+\pi_{\infty,i} -\widetilde \pi_{\infty,ij}=0, \ 1\leq i,j\leq n;\\
			& \overline{A}_{\infty}{\bf 1}^\top\leq {\bf 1}^\top;\\
			& \pi_{\infty,i} \times \Big(\sum_{j=1}^n \overline{a}_{\infty,ij}-1\Big)=0, \ 1\leq i\leq n;\\
			& \overline{A}_{\infty} \ \mbox{is nonnegative};\\
			& \widetilde\pi_{\infty,ij}\overline{a}_{\infty,ij}=0, \ 1\leq i,j\leq n;\\
			&(\pi_{\infty,1},\ldots,\pi_{\infty,n})\in {\cal R}^n_+ \ \mbox{and} \ (\widetilde{\pi}_{\infty,ij})_{n\times n}\in {\cal R}^{n\times n}_+.
		\end{array}
		\right.
		\label{1212-65}
	\end{align}
	To prove (\ref{limit-2}), consider two cases:
	\[
	{\bf Case \ 1}: \ \overline{A}_{\infty}=A_\infty; \ \ {\bf Case \ 2}: \  \overline{A}_{\infty}\neq A_\infty
	\]
	Here $A_\infty$ is given by (\ref{1212-56}), and $\overline{A}_{\infty}$ is given by (\ref{1212-66}).
	
	{\bf Case 1} directly implies (\ref{limit-2}) from (\ref{1212-65}).  Now consider {\bf Case 2}.   From
	the procedure of generating $\{(A_\ell, P_\ell): \ell\geq 0\}$
	and $\{\overline{A}_{\overline{\ell}_k}: k\geq 1\}$, and noting $\{(A_{\overline{\ell}_k},P_{\overline{\ell}_k}): k\geq 0\}\subseteq \{(A_\ell, P_\ell): \ell\geq 0\}$, we have
	\begin{align*}
		\{\overline{A}_{\overline{\ell}_k}: k\geq 1\} \subseteq \{A_\ell: \ell\geq 1\} \ \ \mbox{and } \ \overline{A}_{\overline{\ell}_k}\in
		\{A_{\overline{\ell}_{k-1}+1},\ldots, A_{\overline{\ell}_{k}+1}\}.
	\end{align*}
	This together with (\ref{limit-4}) gives us
	\[
	{\cal L}(A_{\overline{\ell}_{k+1}+1},P_{\overline{\ell}_{k+1}})\leq {\cal L} ( \overline{A}_{\overline{\ell}_{k+1}},P_{\overline{\ell}_{k}})\leq {\cal L}(A_{\overline{\ell}_{k}},P_{\overline{\ell}_{k}}).
	\]
	Letting $k$ approach to infinity, we have ${\cal L}(A_\infty,P_\infty)={\cal L}(\overline{A}_\infty, P_\infty)$.
	This shows us that  ${A}_{\infty}$ solves the convex optimization problem  {\bf (OP[P])} given by (\ref{july-15-6}) with
	$P=P_{\infty}$. Therefore, (\ref{1212-65}) holds with replacing $\overline{A}_\infty$ by $A_\infty$.   (\ref{limit-2})  follows from (\ref{1212-65}).
	Consequently, we have (ii).
	
	Finally  we prove the third part. We consider an arbitrary subsequence $(A_{\ell_k}, P_{\ell_k})$ with limit point $(A_\infty, P_\infty)$. 
	If $( {A}_\infty,{P}_\infty)$ is a maximum point, then there exists $\varepsilon>0$ such that for all feasible $(A,P)$ satisfying $\|A-{A}_\infty\|+\|P-{P}_\infty\|<\varepsilon$, we have ${\cal L}({A}_\infty,{P}_\infty)\geq {\cal L}(A,P)$. Particularly, for feasible $A$ satisfying
	$\|A-{A}_\infty\|<\varepsilon$,
	\begin{align}\label{L_geq}
		{\cal L}({A}_\infty,{P}_\infty)\geq {\cal L}(A,{P}_\infty).
	\end{align}
	
	On the other hand, from the procedure of generating $\{(A_\ell, P_\ell): \ell\geq 0\}$  and $\{(A_{\ell_k},P_{\ell_k}): k\geq 0\}\subseteq \{(A_\ell, P_\ell): \ell\geq 0\}$, we have, for all feasible $A$,
	\[
	{\cal L}(A_{\ell_{k+1}},P_{\ell_{k+1}})\leq {\cal L}(A_{\ell_{k}+1},P_{\ell_{k}})\leq {\cal L}(A,P_{\ell_{k}}).
	\]
	Let $\ell_k$ go to infinity, we obtain
	\begin{align*}
		{\cal L}({A}_\infty,{P}_\infty)\leq {\cal L}(A,{P}_\infty).
	\end{align*}
	This together with \eqref{L_geq} implies that for feasible $A$ satisfying
	$\|A-{A}_\infty\|<\varepsilon$,
	\begin{align}\label{L_equ}
		{\cal L}({A}_\infty,{P}_\infty)= {\cal L}(A,{P}_\infty).
	\end{align}
	Denote such $A$ by ${A}_\infty+(A-{A}_\infty)$ and let $D=A-{A}_\infty$. we can rewrite \eqref{L_equ} as,
	for $D$ satisfying
	\begin{align}
		\left\{
		\begin{array}{lll}
			&\|D\|< \varepsilon;\\
			& {A}_\infty+D \ \mbox{is nonnegative};\\
			& ({A}_\infty+D){\bf 1}^\top \leq {\bf 1}^\top,
		\end{array}
		\right.
		\label{B_con}
	\end{align}
	we have
	\begin{align*}
		{\cal L}({A}_\infty,{P}_\infty)= {\cal L}({A}_\infty+D,{P}_\infty).
	\end{align*}
	In fact, for each $D$ satisfying \eqref{B_con}, $tD$ also satisfies \eqref{B_con} for $t\in [0,1]$. This is because ${A}_\infty+tD=(1-t){A}_\infty+t({A}_\infty+D)$ and ${A}_\infty$ satisfies last two conditions of \eqref{B_con}. This implies that for each $D$ satisfying \eqref{B_con}, we have, for $t\in [0,1]$,
	\begin{align}\label{L_t}
		{\cal L}({A}_\infty,{P}_\infty)= {\cal L}({A}_\infty+tD,{P}_\infty).
	\end{align}
	In the following we prove that \eqref{L_t} is impossible. Recalling \eqref{ob_exp}-\eqref{1212-50},
	\eqref{L_t} is equivalent to
	\begin{align*}
		\sum_{i=1}^n\Big[t^2{P}_{\infty,(i)}DD^\top{P}_{\infty,(i)}^\top+2t{P}_{\infty,(i)}D({A}_\infty^\top-\lambda_iI){P}_{\infty,(i)}^\top \Big]=0.
	\end{align*}
	Since this holds for all $t\in [0,1]$, we must have for $1\leq i\leq n$,
	${P}_{\infty,(i)}DD^\top{P}_{\infty,(i)}^\top=0$ and then ${P}_{\infty,(i)}D={\bf 0}$. Thus we have ${P}_\infty D={\bf 0}$. If we can find a nonsingular $D$ such that \eqref{B_con} holds, then we obtain ${P}_\infty=0$, but this contradicts ${P}_\infty{\bf 1}^\top={\bf 1}^\top.$ This tells us \eqref{L_t} is impossible and we complete the proof. Thus it remains to show there exists nonsingular $D$ such that \eqref{B_con} holds.
	
	Write $A_\infty=(a_{\infty,ij})_{n\times n}$, and let $\delta=\min\{{a}_{\infty,ij}:{a}_{\infty,ij}\neq 0\}\wedge\frac{\varepsilon}{2}\wedge1$ and for $1\leq i,j\leq n$,
	\begin{align*}
		c_{ij}=\left\{
		\begin{array}{rcl}
			\frac{\delta}{2n}, & & \mbox{if} \ {a}_{\infty,ij}=0,\\
			-\frac{\delta}{2}, & & \mbox{if} \ {a}_{\infty,ij}\neq0.
		\end{array}
		\right.
	\end{align*}
	We show that $C=(c_{ij})_{n\times n}$ satisfies
	\begin{align}
		\left\{
		\begin{array}{lll}
			&\|C\|< \frac{\varepsilon}{2};\\
			& {A}_\infty+C \ \mbox{is positive};\\
			& ({A}_\infty+C){\bf 1}^\top < {\bf 1}^\top,
		\end{array}
		\right.
		\label{C_con}
	\end{align}
	The first two conditions are true from the definition of $\delta$. We only need to show the last condition holds. For each $i$, we consider two cases:
	\begin{align*}
		&{\bf Case \ 1^\prime}:   {a}_{\infty, ij}=0 \ \mbox{for all} \ j; \\
		&{\bf Case \ 2^\prime}:   \mbox{at least one of} \ \{ {a}_{\infty,ij},1\leq j\leq n\} \ \mbox{is positive}.
	\end{align*}
	{\bf Case $1^\prime$}: $( {A}_{\infty,(i)}+C_{(i)}){\bf 1}^\top =\sum_{j=1}^nc_{ij}=\frac{\delta}{2}<1$.
	
	\noindent
	{\bf Case $2^\prime$}: $({A}_{\infty,(i)}+C_{(i)}){\bf 1}^\top \leq 1+\sum_{j=1}^nc_{ij}\leq 1-\frac{\delta}{2}+\frac{(n-1)\delta}{2n}<1$.
	
	Thus \eqref{C_con} holds. Then let $D=C+\gamma I$ with $\gamma$ satisfying
	\begin{align}
		\left\{
		\begin{array}{lll}
			&0\leq \gamma<\frac{\varepsilon}{4}\wedge \min_i\{1-(\overline{A}_{(i)}+C_{(i)}){\bf 1}^\top\};\\
			& -\gamma \ \mbox{does not belong to spectrum of $C$}.
		\end{array}
		\right.
		\label{gamma_con}
	\end{align}
	It is easy to verify that $D$ is nonsingular and satisfies \eqref{B_con}. Thus we complete the proof.
	\hfill$\Box$
	
	\section{Numerical Experiments}\label{numerical}
	In this section, we use the {\bf AM} algorithm to conduct the numerical experiments and report the corresponding numerical performance for {\bf Problems 1} and {\bf 2}.
	The numerical experiments have the following implications. One shows that the {\bf AM} algorithm effectively (in sense that the number of the alternating steps ({\bf Steps a} and {\bf b})
	to get a zero-optimal-value is not large) finds the solution for {\bf Problems 1} and {\bf 2} when the problem scale is not large, that is, the {\bf AM} algorithm can quickly produce  $A$ and $P$ such that ${\cal L}(A,P)=0$ and they satisfy the constraints in the nonconvex optimization problem {\bf (OP)} given by (\ref{july-15-2}); and the other indicates that when the {\bf AM} algorithm
	produces a nonzero value, {\bf Problem 1} often has no solution. Thus, the {\bf AM} algorithm provides an efficient approach to simultaneously solve {\bf Problems 1} and {\bf 2} for the small-scale problems. When the problem scale $n$ in {\bf (OP)} given by \eqref{july-15-1} becomes large, due to the number of the variables in
		{\bf (OP)} increases with the rate on a par with $n^2$, it takes a longer time to find the solutions. The following numerical examples show the {\bf AM} algorithm is still implementable when $n$ becomes large around $100$.
	
	
	For all the numerical experiments conducted below, the criterion used to terminate the iterations  is given by
		\begin{align}
			|{\cal L}(A_{k+1},P_{k})-{\cal L}(A_{k},P_{k})|\leq 10^{-6}.\label{alg-term}
		\end{align}
		All the numerical tests are carried out by using python 3.10 on a Dell computer,
		11th Gen Intel(R) Core(TM) i7-11700. We solve the convex optimization problems  {\bf (OP[P])} and {\bf (OP[A])} in each step by using Gurobi Optimizer.

	We first conduct numerical experiments for the case of $n=3$. They demonstrate the initial $A_0$ is critical to the  efficiency of the {\bf AM} algorithm ({\bf Steps a-c}). For different parameter settings, we choose different $A_0$ (denote $\widetilde{A}_0,\overline{A}_0$ and $\widehat{A}_0$) as follows:
	\begin{align*}
		&\widetilde{A}_0=\begin{bmatrix}
			1/3&1/3&1/3\\
			1/3&1/3&1/3\\
			1/6&1/6&1/6
		\end{bmatrix}; \ \ \overline{A}_0=\begin{bmatrix}
			1/9&1/9&1/9\\
			1/9&1/9&1/9\\
			0&1&0
		\end{bmatrix},\\
		& \widehat{A}_0=Q \cdot  \textrm{diag}(\lambda_1,\lambda_2,\lambda_3)\cdot Q^{-1} \ \mbox{and is nonnegative}.
	\end{align*}
	Here $Q$ can be differently chosen while keeping the nonnegativeness of  the matrix $Q\cdot \textrm{diag}(\lambda_1,\lambda_2,\lambda_3)\cdot Q^{-1}.$
	\begin{eqnarray*}
		Q&=&\begin{bmatrix}
			1&1&1\\
			1&1&-1\\
			1&-1&0
		\end{bmatrix} \ \mbox{when $\lambda_1>\lambda_2>\lambda_3\geq 0$};\\
		Q&=&\begin{bmatrix}
			1&1&0\\
			1&-\varepsilon&1-\varepsilon\\
			1&-\varepsilon&-(1-\varepsilon)
		\end{bmatrix} \ \mbox{for sufficiently small $\varepsilon$  when $\lambda_1>\lambda_2\geq0>\lambda_3$};\\
		Q&=&\begin{bmatrix}
			1&1&0\\
			1&-c&1-c\\
			1&-c&-1+c
		\end{bmatrix} \ \mbox{for $c=\frac{\lambda_2}{\lambda_2+\lambda_3}$ when $\lambda_1>0>\lambda_2>\lambda_3$}.
	\end{eqnarray*}
	The idea of choosing $\widetilde{A}_0$ or $\overline{A}_0$ is to make $A_0$ irreducible first, while idea of making $\widehat{A}_0$ to be this form is to directly make its  eigenvalues to be the same as $\{\lambda_i: 1\leq i\leq 3\}$.
	We summarize our numerical test in Table \ref{table-1}. In the table, ``If having an output" means that as long as the algorithm solves equation \eqref{1212-13} exactly while the initial $A_0$ is one of  $\widetilde{A}_0,\overline{A}_0$ and $\widehat{A}_0$, we mark it as $\checkmark$; otherwise, we mark it as \ding{55}.
	Few observations are made as follows:
	\begin{itemize}
		\item For $(\lambda_1,\lambda_2,\lambda_3)=(0.8,-0.2, -0.5)$, and  $(\beta_1,\beta_2,\beta_3)=(0.9375,-0.1042,0.1667)$, if we choose $A_0=\widehat{A}_0$, then the algorithm seems not to produce an optimal zero-value result, however if we choose $A_0=\widetilde{A}_0$, then the algorithm very quickly produces an optimal zero-value result.
		\item Considering $(\lambda_1,\lambda_2,\lambda_3)=(0.8, -0.1,-0.6)$, and $(\beta_1,\beta_2,\beta_3)=(0.9422, 0.0343$, $0.0235)$, if we choose $A_0=\widetilde{A}_0$ or $\widehat{A}_0$, then the algorithm seems not to produce an optimal zero-value result, however if we choose $A_0=\overline{A}_0$, then the algorithm very quickly produces an optimal zero-value result.
		\item Taking
		\begin{equation}\label{examp-nosol}
			\begin{split}
				& (\lambda_1,\lambda_2,\lambda_3)=(0.6,0.4,-0.2), \  (\beta_1,\beta_2,\beta_3)=(0.2727,0.1818,0.5455);\\
				& (\lambda_1,\lambda_2,\lambda_3)=(0.6,0.3,0.1), \  (\beta_1,\beta_2,\beta_3)=(4.7015,-5.3731,1.6716);\\
				& (\lambda_1,\lambda_2,\lambda_3)=(0.8,-0.2, -0.5), \ (\beta_1,\beta_2,\beta_3)=(0.7258, 0.0806, 0.1936),
			\end{split}
		\end{equation}
		regardless of what we choose the initial $A_0$, the {\bf AM} algorithm never produces an optimal zero-value result. Actually, we can also theoretically prove that {\bf Problem 1} has no solution, see a proof in Appendix.
	\end{itemize}
	
	The following example shows that for different initial $A_0$, the {\bf AM} algorithm can give us different {\bf KKT} solutions.
	\begin{myexam}
		Let $(\lambda_1,\lambda_2,\lambda_3)=(0.8,-0.1,-0.6)$, $(\beta_1,\beta_2,\beta_3)=(0.9422,0.0343,0.0235)$,
		and $A_0$ be one of 
		\begin{align*}
			\widetilde{A}_0=\begin{bmatrix}
				1/3&1/3&1/3\\
				1/3&1/3&1/3\\
				1/6&1/6&1/6
			\end{bmatrix}, \
			 \widehat{A}_0=\begin{bmatrix}
				0.0125&0.39375&0.39375\\
				0.1125&0.04375&0.64375\\
				0.1125&0.64375&0.04375
			\end{bmatrix}, \
			\overline{A}_0=\begin{bmatrix}
				1/9&1/9&1/9\\
				1/9&1/9&1/9\\
				0&1&0
			\end{bmatrix}.
		\end{align*} 
		For initial $A_0=\widetilde{A}_0$, the outcome of the {\bf AM} algorithm is
		{\rm \begin{align*}
			& \widetilde{P}_\infty \approx \begin{bmatrix}
				1.5931\text{e-}04&9.9968\text{e-}01&1.5931\text{e-}04\\
				4.9686\text{e-}01&6.2741\text{e-}03&4.9686\text{e-}01\\
				4.5738\text{e-}01&8.5234\text{e-}02&4.5738\text{e-}01
			\end{bmatrix}, \
			\widetilde{A}_\infty \approx \begin{bmatrix}
				0&0&0\\
				0&7.8986\text{e-}01&0\\
				0&0&0
			\end{bmatrix},\\
			&{\cal L}(\widetilde{A}_\infty,\widetilde{P}_\infty) \approx 1.6973\text{e-}01; 
		\end{align*} }
		for initial $A_0=\widehat{A}_0$, the outcome is
		{\rm  \begin{align*}
			&\widehat{P}_\infty \approx \begin{bmatrix}
				1.6051\text{e-}04&1.6051\text{e-}04&9.9968\text{e-}01\\
				4.9686\text{e-}01&4.9686\text{e-}01&6.2747\text{e-}03\\
				4.5738\text{e-}01&4.5738\text{e-}01&8.5239\text{e-}02
			\end{bmatrix},
			\widehat{A}_\infty \approx \begin{bmatrix}
				0&0&0\\
				0&0&0\\
				0&0&7.8986\text{e-}01
			\end{bmatrix},\\
			&{\cal L}( \widehat{A}_\infty, \widehat{P}_\infty) \approx 1.6973\text{e-}01;
		\end{align*}}
		for initial $A_0=\overline{A}_0$, the outcome is
	{\rm 	\begin{align*}
			&\overline{P}_\infty \approx \begin{bmatrix}
				1.6057\text{e-}01&4.3229\text{e-}01&4.0714\text{e-}01\\
				-1.0098 &2.2643 &-2.5447\text{e-}01\\
				6.3827\text{e-}02&3.2933 &-2.3571 
			\end{bmatrix},\\
			&\overline{A}_\infty \approx \begin{bmatrix}
				8.8886\text{e-}02&1.2552\text{e-}02&8.9856\text{e-}01\\
				1.0247\text{e-}01&5.5737\text{e-}09&4.1210\text{e-}01\\
				1.6164\text{e-}01&8.3836\text{e-}01&2.3440\text{e-}07
			\end{bmatrix},\\
			&{\cal L}(\overline{A}_\infty,\overline{P}_\infty) \approx 3.4547\text{e-}05.
		\end{align*} }
		This numerical example shows us that if we choose $A_0=\overline{A}_0$, we can obtain a global minimum $(\overline{A}_\infty,\overline{P}_\infty);$ but if we choose $A_0=\widetilde{A}_0$ or $A_0=\widehat{A}_0$, we can only obtain a {\bf KKT} solution. \hfill$\Box$
	\end{myexam}

		The following is a concrete example to further illustrate Remark \ref{rem-sequence-convergence}.
		\begin{myexam}\label{examp-noconverge}
			Let $(\lambda_1,\lambda_2,\lambda_3)=(0.8,0.3637, 0.5120)$, and $(\beta_1,\beta_2,\beta_3)=(0.1700, 0.1500$, $0.6800)$.
			Starting with $A_0={\bf 0}$ to initiate {\bf Step-a} for {\bf (OP[A])}, we have
			\begin{align*}
				& P_0 \approx \begin{bmatrix}
					0.3333&0.3333&0.3333\\
					0.3333&0.3333&0.3333\\
					0.3333&0.3333&0.3333
				\end{bmatrix}.
			\end{align*}
			Then we move to {\bf Step-b} for {\bf (OP[P])} with this $P_0$. We obtain two solutions$:$
			{\rm  \begin{align*}
				&A_1 \approx \begin{bmatrix}
					1.8618\text{e-}01&1.8618\text{e-}01&1.8618\text{e-}01\\
					1.8618\text{e-}01&1.8618\text{e-}01&1.8618\text{e-}01\\
					1.8618\text{e-}01&1.8618\text{e-}01&1.8618\text{e-}01
				\end{bmatrix}, \\
				&\widehat{A}_1 \approx \begin{bmatrix}
					1.8618\text{e-}01&1.8618\text{e-}01&1.8618\text{e-}01\\
					1.8618\text{e-}01&1.8618\text{e-}01&3.7219\text{e-}01\\
					1.8618\text{e-}01&1.8618\text{e-}01&1.8889\text{e-}04
				\end{bmatrix}.
			\end{align*}} 
   Moving forward to {\bf Step-c} to continue the {\bf AM} algorithm from $A_1$, it turns out that $P_\ell \approx P_0$ and $A_{\ell+1} \approx A_1$ for all $\ell\geq 1$, this implies that
			$\{(A_\ell,P_\ell): \ell\geq 0\}$ has a limit point $(A_\infty,P_\infty) \approx (A_1,P_0)$.
			
			Next we move to {\bf Step-c} to continue the {\bf AM} algorithm but from $\widehat{A}_1$. We generate another sequence
			$\{(A_0,P_0), (\widehat{A}_\ell,\widehat{P}_\ell): \ell\geq 1\}$. It admits a limit point$:$
		{\rm  \begin{align*}
				&\widehat{P}_\infty \approx \begin{bmatrix}
					3.2331\text{e-}01&4.6428\text{e-}01&2.1241\text{e-}01\\
					5.6749\text{e-}02&-5.2850\text{e-}02&9.9610\text{e-}01\\
					4.0872\text{e-}01&-5.2398\text{e-}02&6.4368\text{e-}01
				\end{bmatrix},\\
				&\widehat{A}_\infty \approx \begin{bmatrix}
					5.3558\text{e-}01&5.1538\text{e-}05&2.5786\text{e-}01\\
					1.8413\text{e-}01&7.8962\text{e-}01&2.6125\text{e-}02\\
					7.3492\text{e-}06&2.2587\text{e-}02&3.5039\text{e-}01
				\end{bmatrix}.
		\end{align*}} 
			This implies that starting from the same $A_0$, due to the nonuniqueness of solutions to {\bf (OP[P])}, the {\bf AM} algorithm can generate different sequences whose limit points are also different. \hfill$\Box$
		\end{myexam}

	\begin{table*}[htb]
		\center
		\caption{Numerical Experiments}
		\begin{threeparttable}
			\renewcommand{\arraystretch}{2.6}
			\resizebox{\linewidth}{!}{
				\begin{tabular}{|c|c|c|c|c|c|c|c|c|} \hline
					$\lambda_1$ &	$\lambda_2$	&$\lambda_3$& $\beta_1$	&$\beta_2$&	$\beta_3$& If having an output & Performance of $\widetilde{A}_0$& Performance of $\widehat{A}_0$\\ \hline
					
					0.6 &0.4&	-0.3	&0.3786	&0.5049&	0.1165&	 $\checkmark$	&55 steps&	2038 steps\\ \hline
					
					0.6	&0.4&-0.2&	0.5&	0.3333&	0.1667 &		$\checkmark$	&85 steps&	stuck, cannot decrease to 0\\ \hline

					0.6&	0.3	&0.1&	0.4960	&0.2835&	0.2205 &		$\checkmark$	&19 steps&	33 steps\\ \hline
					
					0.6&	0.3&	0.1&	0.4455	&0.3564&	0.1981  &		$\checkmark$	&24 steps&	30 steps\\ \hline
					
					0.6&	0.4&-0.2&	0.2727&	0.1818&	0.5455  &		\ding{55}	&stuck, cannot decrease to 0&	stuck, cannot decrease to 0\\ \hline
					
					0.6	&0.3&	0.1&	4.7015&	-5.3731	&1.6716 &		\ding{55}	&stuck, cannot decrease to 0&	stuck, cannot decrease to 0\\ \hline

					0.8 &-0.2&	-0.5	&0.9375	&-0.1042&	0.1667&	 $\checkmark$	&190 steps&	stuck, cannot decrease to 0\\ \hline
					
					0.8	&-0.2&-0.5&	0.7759&	0.0862&	0.1379 &		$\checkmark$	&199 steps&	stuck, cannot decrease to 0\\ \hline

					0.8&	-0.3&-	0.4&	0.9440	&0.0290&	0.0270  &		$\checkmark$	&298 steps&	stuck, cannot decrease to 0\\ \hline
					
					0.8	&0&	-0.4&	0.9358&	0.0375	&0.0267 &		$\checkmark$	&64 steps&stuck, cannot decrease to 0\\ \hline
					
					0.8&	-0.1	&-0.6&	0.9422	&0.0343&	0.0235 &	{\color{red}$\checkmark $}$^1$	&stuck, cannot decrease to 0&	stuck, cannot decrease to 0\\ \hline

					0.8&	-0.2&-0.5&	0.7258&	0.0806&	0.1936  &\ding{55}	&stuck, cannot decrease to 0&	stuck, cannot decrease to 0\\ \hline

					0.5 &0.3&	-0.4	&1.6716	&-0.2985&	-0.3731&	 $\checkmark$	&stuck, cannot decrease to 0&	1169 steps\\ \hline
					
					0.5	&0.3&-0.4&	1.1546&-0.2062&	0.0516 &		$\checkmark$	&22 steps&	202 steps\\ \hline
					
					0.5&	0.3&-0.4&	0.8175&	0.1460&	0.0365  &		$\checkmark$	&144 steps &	214 steps\\ \hline

					0.5&	0.1	&-0.3&	0.5969	&0.2653&	0.1378 &		$\checkmark$	&72 steps&	39 steps\\ \hline
					
					0.5&	0.1&	-0.3&	0.1297	&0.7206&	0.1497  & 	$\checkmark$	&11427 steps&	stuck, cannot decrease to 0\\ \hline
					
					0.5	&0.1&	-0.2&	0.0496&	0.8264	&0.1240 &	 	$\checkmark$	&1785 steps&	stuck, cannot decrease to 0\\ \hline
				\end{tabular}
			}	
			\begin{tablenotes}
				\footnotesize
				\item[1.] If we choose initial $A_0=\overline{A}_0$, then we have an output.
			\end{tablenotes}
		\end{threeparttable} \label{table-1}
	\end{table*}

		The following numerical example demonstrates that the minimal phase-type representation problem discussed in Remark \ref{rem-1212-2} can be solved by {\bf Problem 1}.	
		\begin{myexam}
			Given a rational function
			\begin{align*}
				f(z)=z \Big(\frac{0.2927}{1-0.8000z}-\frac{0.0376}{1-0.7095z}-\frac{0.0523}{1-0.3473z}-\frac{0.0673}{1-0.3246z}-\frac{0.1213}{1-0.2132z}\Big),
			\end{align*}
			we want to determine whether there exists a discrete-time Markov chain with state space $\{0,1,\ldots,5\}$ such that
			{\rm (i)} state $0$ is an absorbing state and the other states are transient; and  {\rm (ii)} the moment generating function of its absorbing time is $f(z)$.
			
			Note that $f(z)$ has the same form as \eqref{1212-14}. Then according to \eqref{beta} and Remark {\rm \ref{rem-1212-2}},
			letting
			\begin{align*}
				&\lambda=(0.8000,0.7095,0.3473,0.3246,0.2132),\\
				&\beta=(1.4635,-0.1295,-0.0802,-0.0996,-0.1542),
			\end{align*}
			this problem is equivalent to solve {\bf Problem 1} with these $\lambda$ and $\beta$.
			Now we use the {\bf AM} algorithm with beginning {\rm $A_0= \textrm{diag}(\lambda)$} to solve the optimization problem {\bf (OP)} given by \eqref{july-15-2}.	
			The outcome of the {\bf AM} algorithm is
		{\rm 	\begin{align*}
				&P_{\infty} \approx \begin{bmatrix}
					6.7303\text{e-}01  &  8.8669\text{e-}02 &   6.6836\text{e-}02  &  6.7412\text{e-}02  &  1.0405\text{e-}01 \\
					-2.0248\text{e-}03  &  1.0021 &   4.3047\text{e-}04 &  1.0153\text{e-}05  &  -5.4163\text{e-}04 \\
					-1.6450\text{e-}03  &  -3.3601\text{e-}04 &   1.0161 & -1.1237\text{e-}02  &  -2.8517\text{e-}03 \\
					-1.8440\text{e-}03 &   -1.3702\text{e-}04 &  1.3211\text{e-}02  &  9.9301\text{e-}01  & -4.2388\text{e-}03 \\
					4.4448\text{e-}06  &  2.1262\text{e-}04 &  3.3519\text{e-}03  &  4.2398\text{e-}03  &  9.9219\text{e-}01
				\end{bmatrix},\\
				&A_{\infty} \approx \begin{bmatrix}
					7.9981\text{e-}01 &   1.1896\text{e-}02 &   4.4917\text{e-}02  &  4.7594\text{e-}02  &  9.0717\text{e-}02  \\
					1.8482\text{e-}04  & 7.0952\text{e-}01  &  2.4698\text{e-}04   & 1.1118\text{e-}04   & 1.4087\text{e-}05   \\
					7.4225\text{e-}04  &  1.4375\text{e-}04  & 3.4737\text{e-}01   & 2.1317\text{e-}06   & 6.2076\text{e-}07   \\
					8.7259\text{e-}04  &  8.7464\text{e-}05  &  8.7411\text{e-}07  &3.2469\text{e-}01   & 1.6821\text{e-}07   \\
					3.4582\text{e-}05  &  7.5530\text{e-}06  &  3.6994\text{e-}08   & 3.4510\text{e-}08   &2.1320\text{e-}01   \\
				\end{bmatrix},
			\end{align*}}
			and {\rm ${\cal L}(A_{\infty},P_{\infty}) \approx  6.7854\text{e-}07$}. This implies that {\bf Problem 1} has a solution. From Remark {\rm \ref{rem-1212-2}}, let {\rm $$\alpha=\beta P_{\infty} \approx  (9.8556\text{e-}01,7.0066\text{e-}09,1.4438\text{e-}02,3.1963\text{e-}08,1.9817\text{e-}10).$$}	
			We obtain a discrete-time Markov chain with initial distribution $\alpha$ and transition probability matrix $A_{\infty}$ over transient states $\{1,\ldots,5\}$.
			\hfill$\Box$	
		\end{myexam}

		To show the {\bf AM} algorithm can also solve the larger problems, now we conduct numerical experiments by selecting $n=20, 40, 60, 80$ and $100$.
		For each fixed $n$, we do the following:
		\begin{itemize}
			\item [(i)] randomly select parameters $\lambda$ and $\beta$ but keeping \eqref{1212-11}-\eqref{1212-12} hold;
			\item [(ii)] terminate the alternating iterations from the {\bf AM} algorithm by the criterion given in \eqref{alg-term}, say the termination happens at alternating iteration step $k$;
			\item [(iii)] calculate the corresponding objective function value ${\cal L}(A_{k},P_{k})$;
			\item[(iv)] classify ${\cal L}(A_{k},P_{k})$ calculated in (iii) into two groups, if ${\cal L}(A_{k},P_{k})\leq10^{-4}$, put it to group-1 ({\bf Problem 1} has a solution), otherwise, put it to group-2 ({\bf Problem 1} may not have a solution);
			\item[(v)] repeat (i)-(iv) until each group has 20 ${\cal L}(A_{k},P_{k})$'s with the convention that if one group, say group-1, first receives 20 ${\cal L}(A_{k},P_{k})$'s before group-2 receives 20 ${\cal L}(A_{k},P_{k})$'s, then we keep doing (i)-(iv) to find ${\cal L}(A_{k},P_{k})$ larger than $10^{-4}$ until group-2 gets 20 ${\cal L}(A_{k},P_{k})$'s while we ignore ${\cal L}(A_{k},P_{k})$'s that are less than or equal to $10^{-4}$ and stop adding them to group-1.
		\end{itemize}
		First we start with $A_0= \textrm{diag}(\lambda)$, group-1 and group-2 are reported in Tables \ref{table:1} and \ref{table:3}, respectively. The data includes the average computing time, the average number of iterations, and the average objective function value from the 20 repetitions. When $n=80, 100$, we cannot get 20 ${\cal L}(A_{k},P_{k})$'s which are larger than $10^{-4}$ as they take more than 8 hours, and still cannot get convergent values.
		\begin{table}[htb]
			\begin{center}
				\caption{Numerical experiments when {\bf Problem 1} has a solution}
				\label{table:1}
				\begin{tabular}{|m{1.5cm}<{\centering}|m{2.6cm}<{\centering}|m{4cm}<{\centering}|m{5cm}<{\centering}|}
					\hline   ${\bf n}$ & \textbf{Computing time} & \textbf{Number of the iterations} &\textbf{Final objective function value}\\
					\hline   20 & 0.9192 s & 2.95  &5.8333e-07 \\
					\hline   40  & 10.5660 s & 2.3  &8.3175e-07\\
					\hline   60 & 44.9387 s & 2.05  &6.1865e-07\\
					\hline   80 & 181.3309 s & 2.0  & 4.3424e-07 \\
					\hline   100 & 1589.4458 s & 2.0  &2.5225e-07\\
					\hline
				\end{tabular}
			\end{center}
		\end{table}
		\begin{table}[htb]
			\begin{center}
				\caption{Numerical experiments for $A_0= \textrm{diag}(\lambda)$}
				\label{table:3}
				\begin{tabular}{|m{1.5cm}<{\centering}|m{2.6cm}<{\centering}|m{4cm}<{\centering}|m{5cm}<{\centering}|}
					\hline   ${\bf n}$ & \textbf{Computing time} & \textbf{Number of the iterations} &\textbf{Final objective function value}\\
					\hline   20 & 43.2514 s & 147.9  &0.0878 \\
					\hline   40  & 531.7167 s & 113.4  &0.1105\\
					\hline   60 & 2159.3460 s & 95.0  &0.1138\\
					\hline   80 & $> 8$ h & $> 100$  & - \\
					\hline   100 & $> 8$ h & $> 100$  & - \\
					\hline
				\end{tabular}
			\end{center}\label{sln-no-exist}
		\end{table}
		
		Next we start with $A_0={\bf 0}$ to generate another group-1 and group-2. It turns out that group-1 has a similar result with the previous group-1 ($A_0= \textrm{diag}(\lambda)$, Table
		\ref{table:1}) while group-2, unlike the previous group-2 (Table \ref{table:3}), can be quickly generated with 20 ${\cal L}(A_{k},P_{k})$'s, we report them in Table \ref{table:2}.
		
		
		
		
		From these numerical analyses, we observe the following:
		\begin{itemize}
			\item Table \ref{table:1} shows that when the parameters $\lambda$ and $\beta$ make {\bf Problem 1} have a solution, the algorithm requires few iterations, and quickly provides a solution;
			\item Tables \ref{table:3} and \ref{table:2} show that when the parameters $\lambda$ and $\beta$ may make {\bf Problem 1} have no solution, the efficiency of the {\bf AM} algorithm is affected much by the initial $A_0$ selection. This is consistent with what we observe in Table \ref{table-1};
			\item The computing time increases rapidly as $n$ increases regardless of whether {\bf Problem 1} has a solution. This is because the number of variables in {\bf (OP)} is $2n^2$.
		\end{itemize}

		\begin{table}[htb]
			\begin{center}
				\caption{Numerical experiments for $A_0={\bf 0}$}
				\label{table:2}
				\begin{tabular}{|m{1.5cm}<{\centering}|m{2.6cm}<{\centering}|m{4cm}<{\centering}|m{5cm}<{\centering}|}
					\hline   ${\bf n}$ & \textbf{Computing time} & \textbf{Number of the iterations} &\textbf{Final objective function value}\\
					\hline   20 & 1.5537 s & 4.3  &0.2109 \\
					\hline   40  & 8.0633 s & 1.6  &0.2217\\
					\hline   60 & 31.2986 s & 1.4  &0.2258\\
					\hline   80 & 85.7927 s & 1.4  & 0.2256 \\
					\hline   100 & 909.4467 s & 1.5  &0.2256\\
					\hline
				\end{tabular}
			\end{center}\label{sln-no-exist-1}
		\end{table}

	\section{Concluding Remarks}\label{conc}
	In this paper we investigate the substochastic inverse eigenvalue problem (SstIEP) with the corresponding eigenvector constraints by constructing a substochastic matrix.
	The eigenvector constraints make the constructed substochastic matrix to be the transition probability matrix of a finite-state discrete-time Markov chain whose each state is transient except one absorbing state. This SstIEP is equivalently transformed into a nonconvex optimization problem (NcOP).
	The equivalence means that when NcOP has a zero optimal value, the SstIEP is solvable, and the solution of NcOP is what we want to find the substochastic matrix; and when NcOP has a nonzero optimal value, the SstIEP is insolvable.  The  {\bf AM} algorithm is developed for NcOP to be used constructing the desired substochastic matrix. The proposed algorithm is proved to converge for any initial value, and numerically shown to be quite effective. At the same time, the numerical experiments show us the initial value is sometimes critical to find a zero optimal value. So the next research is to look at how to optimally choose the initial value in order to more efficiently find a zero optimal value, and consider the case when the spectrum data $\{\lambda_1,\lambda_2,\ldots,\lambda_n\}$ is a general self-conjugate complex numbers. The latter case is difficult because (i) if some $\lambda_i$ is multiple (its multiplicity is more than one), then $ \textrm{diag}(\lambda)$ in \eqref{1212-13} is changed into the Jordan canonical form of $A$. Since the Jordan form with respect to the multiple eigenvalue is not unique, this problem cannot be formulated as a clean form \eqref{1212-13} without a specific Jordan form; and (ii) the proof of Proposition \ref{tech-prop} heavily depends on the structure of $B^\top$ given by \eqref{12-10-1}. If the Jordan form is not diagonal, we cannot obtain the result similar to Proposition \ref{tech-prop} using the same method, and then we cannot obtain the convergence of the proposed algorithm.

	\section*{Acknowledgments}
	The authors thank Associate Editor, Prof. Beatrice Meini, and two anonymous referees for their numerous thoughtful and
		constructive comments which have helped them to significantly improve the paper.
	The authors are grateful to Yinyu Ye at Stanford University and Chao Ding at Academy of Mathematics and Systems Science, Chinese Academy of Sciences,
	for many fruitful discussions on the KKT  conditions and their convergence.

	\newpage
	
		\section{Appendix}
		\renewcommand{\theequation}{A-\arabic{equation}}
		\setcounter{equation}{0}
		We show that\\[0.15in]
		\noindent
		{\bf Claim A.} {\bf Problem 1} {\it  does not have a solution for $\lambda$ and $\beta$ given by \eqref{examp-nosol}}.\\[0.15in]
		Here we consider a general $n$. We first establish that {\bf Problem 1} is equivalent to a class of the positive realization problems in control theory.
		Then the equivalence makes us to directly apply the extant results in the positive realization problems to {\bf Problem 1}. For
		\begin{align}\label{hz}
			H(z)=\frac{r_1}{z-\lambda_1}+\frac{r_2}{z-\lambda_2}+\cdots+\frac{r_n}{z-\lambda_n} \ \mbox{with $r_i\neq 0$ for $1\leq i\leq n$},
		\end{align}
		we say $H(z)$ has an $n$-dimensional positive realization $(A,b,c)$ if there exist $(A,b,c)$ with $A \in {\cal R}_{+}^{n\times n}$, $b,c \in {\cal R}_{+}^{n}$ such that $H(z)=c(zI-A)^{-1}b^{\top}$.
		
		The following proposition establishes the equivalence between the solution existence of {\bf Problem 1} and the positive realization problem, and its proof is based on the two equivalences: one is between {\bf Problem 1} and the minimal phase-type representation, which is given by Remark \ref{rem-1212-2}; and the other one is about the minimal phase-type representation and the positive realization problem, which is given by \cite{liu2023}.
		The detailed proof is omitted here. \\[0.15in]
		\noindent
		{\bf Proposition A1.}	{\it Let $r_i=(1-\lambda_i)\beta_i$ for $1\leq i\leq n$. Then {\bf Problem 1} with $\lambda$ and $\beta$ has a solution if and only if $H(z)$ given by \eqref{hz} has an $n$-dimensional positive realization.}\\[0.15in]
		With Proposition A1, we can now proceed to prove {\bf Claim A}.\\[0.1in]
		\noindent
		{\it Proof of \ {\bf Claim A:}} First by Proposition A1, we only need to show that there does not exist 3-dimensional positive realization for $H(z)$ given by \eqref{hz} with $n=3$ when $\lambda$ and $\beta$ satisfy \eqref{examp-nosol}.  The nonexistence of the 3-dimensional positive realization for such $\lambda$ and $\beta$ can be directly proven using Theorems 1 and 2 of \cite{benvenuti2012minimal} by verifying  that
		\[(\lambda_1,\lambda_2,\lambda_3)=(0.6,0.4,-0.2) \ \mbox{and  } (\beta_1,\beta_2,\beta_3)=(0.2727,0.1818,0.5455)\]
		make $(1-\lambda_1)\beta_1\lambda_1+(1-\lambda_2)\beta_2\lambda_2+(1-\lambda_3)\beta_3\lambda_3<0$;
		\[
		(\lambda_1,\lambda_2,\lambda_3)=(0.6,0.3,0.1)  \ \mbox{and  }  (\beta_1,\beta_2,\beta_3)=(4.7015,-5.3731,1.6716)
		\]
		have $(1-\lambda_1)\beta_1+(1-\lambda_2)\beta_2+(1-\lambda_3)\beta_3<0$, and
		\[ (\lambda_1,\lambda_2,\lambda_3)=(0.8,-0.2, -0.5)  \ \mbox{and  }  (\beta_1,\beta_2,\beta_3)=(0.7258, 0.0806, 0.1936)
		\]
		enable $(1-\lambda_1)\beta_1\lambda_1+(1-\lambda_2)\beta_2\lambda_2+(1-\lambda_3)\beta_3\lambda_3<0$.
		\hfill $\Box$
		
	\end{document}